\documentclass[11pt]{article}
\usepackage[english]{babel}
\usepackage{amsmath,amsthm,amssymb}
\usepackage{mathrsfs}
\usepackage{bbm}
\usepackage{amsfonts}
\usepackage[margin=0.8in]{geometry}
\usepackage{graphicx}
\usepackage{caption}
\usepackage{subcaption}
\usepackage{wasysym}
\usepackage{enumerate}
\usepackage{algorithm2e}
\usepackage{tabularx}
\usepackage{color}
\usepackage{wrapfig}
\usepackage{todonotes}

\newtheorem{thm}{Theorem}[section]
\newtheorem{ass}[thm]{Assumption}
\newtheorem{cor}[thm]{Corollary}
\newtheorem{lem}[thm]{Lemma}
\newtheorem{prop}[thm]{Proposition}
\newtheorem*{hyp*}{Hypothesis}
\theoremstyle{definition}
\newtheorem{defn}[thm]{Definition}

\theoremstyle{rem}
\newtheorem{rem}[thm]{Remark}
\newtheorem{example}[thm]{Example}
\numberwithin{equation}{section}

\newcommand{\R}{\mathbb R}

\newcommand{\bbD}{\mathbb D}
\newcommand{\bbF}{\mathbb F}

\newcommand{\bbG}{\mathbb G}

\newcommand{\mcA}{\mathcal{A}}

\newcommand{\mcB}{\mathcal{B}}
\newcommand{\mcD}{\mathcal D}
\newcommand{\mcE}{\mathcal E}
\newcommand{\mcG}{\mathcal{G}}
\newcommand{\mcN}{\mathcal N}
\newcommand{\mcF}{\mathcal F}

\newcommand{\mcT}{\mathcal T}
\newcommand{\mcP}{\mathcal P}

\newcommand{\mcH}{\mathcal H}

\newcommand{\mcU}{\mathcal U}
\newcommand{\mcS}{\mathcal S}
\newcommand{\Sqlc}{\mathcal{S}_{\textit{q}}}

\newcommand{\E}{\mathbb{E}}
\newcommand{\Prob}{\mathbb{P}}

\newcommand{\vect}{\bold{t}}
\newcommand{\vecb}{\bold{b}}

\newcommand{\vecv}{\bold{v}}
\newcommand{\esssup}{\mathop{\rm{ess}\,\sup}}

\newcommand{\argmin}{\mathop{\arg\min}}
\newcommand{\ett}{\mathbbm{1}}

\newcommand{\cadlag}{c\`adl\`ag~}
\newcommand{\cadlagSTOP}{c\`adl\`ag}

\newcommand{\cadlagp}{c\`adl\`ag.~}

\newcommand{\ie}{\textit{i.e.\ }}
\newcommand{\eg}{\textit{e.g.\ }}

\begin{document}

\title{Infinite Horizon Impulse Control of Stochastic Functional Differential Equations\footnote{This work was supported by the Swedish Energy Agency through grant number 48405-1}}

\author{Magnus Perninge\footnote{M.\ Perninge is with the Department of Physics and Electrical Engineering, Linnaeus University, V\"axj\"o,
Sweden. e-mail: magnus.perninge@lnu.se.}} %
\maketitle
\begin{abstract}
We consider impulse control of stochastic functional differential equations (SFDEs) driven by L\'evy processes under an additional $L^p$-Lipschitz condition on the coefficients. Our results, which are first derived for a general stochastic optimization problem over infinite horizon impulse controls and then applied to the case of a controlled SFDE, apply to the infinite horizon as well as the random horizon settings. The methodology employed to show existence of optimal controls is a probabilistic one based on the concept of Snell envelopes.
\end{abstract}

\section{Introduction}
The standard stochastic impulse control problem is an optimal control problem that arises when an operator controls a dynamical system by intervening on the system at a discrete set of stopping times. Generally, an intervention can be represented by an element in the control set $U$ which we assume to be a compact subset of $\R^m$.

In impulse control the control-law, thus, takes the form $u=(\tau_1,\ldots,\tau_N;\beta_1,\ldots,\beta_N)$, where $\tau_1\leq\tau_2\leq\cdots\leq\tau_N$ is a sequence of times when the operator intervenes on the system and $\beta_j$ is the impulse that the operator affects the system with at time $\tau_j$. The standard impulse control problem in infinite horizon can be formulated as finding a control that maximizes
\begin{equation}\label{ekv:std-imp-cont}
\E\bigg[\int_0^\infty e^{-\rho_0 s}\phi(s,X^u_s)ds-\sum_{j=1}^Ne^{-\rho_0 \tau_j}\ell(\tau_j,X^u_{\tau_j},\beta_j)\bigg],
\end{equation}
where $\rho_0>0$ is a constant referred to as the discount factor, $X^u$ is an $\R^m$-valued controlled stochastic process that jumps at interventions (\eg by setting $X^u_{\tau_j}=\Gamma(\tau_j,X^u_{\tau_j-},\beta_j)$ for some deterministic function $\Gamma$) and the deterministic functions\footnote{Throughout, we let $\R_+:=[0,\infty)$} $\phi:\R_+\times\R^m\to\R$ and $\ell:\R_+\times \R^m\times U\to\R$ give the running reward and the intervention costs, respectively. The quantity $\ell(t,x,b)$, thus, represents the cost incurred by applying the impulse $b\in U$ at time $t\in \R_+$ when the state is $x$.

As impulse control problems appear in a vast number of real-world applications (see \eg~\cite{Korn,PSImpulsive} for applications in finance and \cite{BaseiImpulse,CarmLud} for applications in energy) a lot of attention has been given to various types of problems where the control is of impulse type. In the standard Markovian setting, where $X^u$ solves a stochastic differential equation (SDE) driven by a L\'evy process on $[\tau_j,\tau_{j+1})$, the relation to quasi-variational inequalities has frequently been exploited to find optimal controls (see the seminal work in~\cite{BensLionsImpulse} or turn to \cite{OksenSulemBok} for a more recent textbook). In the non-Markovian framework an impulse control problem in finite horizon ($T<\infty$) was solved in~\cite{DjehiceImpulse} by utilizing the link between optimal stopping and reflected BSDEs (originally discovered in~\cite{ElKaroui1}) while considering the reward functional
\begin{equation*}
\E\bigg[\int_0^T \phi(s,\omega,L_t^u)ds-\sum_{j=1}^Nc(\beta_j)\bigg],
\end{equation*}
where $\phi:[0,T]\times\Omega\times \R^n\to\R$ is now a random (and not necessarily Markovian) field and the controlled process $L^u$ takes the particular form $L^u_t:=L_t+\sum_{j= 1}^N\ett_{[\tau_j\leq t]}\beta_j$, with $L$ an (exogenous) non-controlled process and assuming that $U$ is a finite set. Relevant is also the treatment of multi-modes optimal switching problems in a non-Markovian setting in~\cite{BollanMSwitch1}.

Almost all production systems are subject to delays in the sense that some time is required to start up the production units. In this regard a lot of effort has been directed at impulse controls where the effect of the interventions are delayed by a fixed lag. In the Markovian setting, the novel paper~\cite{BarIlan95} proposes an explicit solution to an inventory problem with uniform delivery lag by taking the current stock plus pending orders as one of the states. Similar approaches are taken in \cite{RAid2} where explicit optimal solutions of impulse control problems with uniform delivery lags are derived for a large set of different problems and~\cite{Bruder} that propose an iterative algorithm. The authors of \cite{OksenImpulse} propose a solution to an impulse control problem with uniform delivery lags when interventions are prohibited during the delay period of prior interventions. Their approach is based on defining an operator that circumvents the delay period.

Also in the non-Markovian setting problems with delivery lag have been considered. In~\cite{Hdhiri} the original work of \cite{DjehiceImpulse} was extended to incorporate delivery lag by setting  $L^u_t:=L_t+\sum_{j= 1}^N\ett_{[\tau_j+\Delta\leq t]}\beta_j$ for a fixed $\Delta>0$. As in \cite{OksenImpulse} the work in \cite{Hdhiri} is based on the assumption that $\tau_{j+1}\geq\tau_j+\Delta$. This work was later extended by considering the infinite horizon setting in~\cite{DjehicheInfHorImp}. Notable is also the recent work on finite horizon impulse control of SFDEs driven by a Brownian motion in~\cite{JonteSFDE}.

In the present article we take a different approach to all the above mentioned works by considering the abstract reward functional
\begin{align}
J(u)&:=\E\bigg[\varphi(\tau_1,\ldots,\tau_{N};\beta_1,\ldots,\beta_{N}) -\sum_{j=1}^N c(\tau_1,\ldots,\tau_{j};\beta_1,\ldots,\beta_{j})\bigg],\label{ekv:objfun}
\end{align}
where the terminal reward $\varphi$ maps controls to values of the real line and is measurable with respect to $\mcG\otimes\mcB(\mcD)$, where $\mcB(\mcD)$ is the Borel $\sigma$-field of $\mcD:=\cup_{i=0}^\infty D^i$, with $D^i:=\{(t_1,\ldots,t_i;b_1,\ldots,b_i): 0\leq t_1\leq\cdots\leq t_i,\, b_j\in U\}$ and $\mcG$ is the $\sigma$-field of a complete probability space $(\Omega,\mcG,\Prob)$. The intervention cost $c$ is also assumed to be a $\mcG\otimes\mcB(\mcD)$-measurable map in addition to being bounded from below by a deterministic positive function. We consider the partial information setting and assume that we observe the system through a filtration $\bbF:=\{\mcF_t\}_{t\geq 0}$ of sub-$\sigma$-fields of $\mcG$ and thus restrict our attention to $\bbF$-adapted controls.

To indicate the applicability of the results we consider the special case when
\begin{align}\label{ekv:phiSFDE}
\varphi(u)&=\int_0^\infty e^{-\rho(s)}\phi(s,X^u_s)ds,
\end{align}
and
\begin{align}\label{ekv:cSFDE}
c(\tau_1,\ldots,\tau_{j};\beta_1,\ldots,\beta_{j})&= e^{-\rho(\tau_j)}\ell(\tau_j,X^{(\tau_1,\ldots,\tau_{j-1};\beta_1,\ldots,\beta_{j-1})}_{\tau_j},\beta_j),
\end{align}
where $X^u$ solves an impulsively controlled stochastic functional differential equation (SFDE) driven by a L\'evy process under an additional $L^p$-type Lipschitz condition on the coefficients of the SFDE. Furthermore, we will see that the results easily extend to problems with a random horizon which allows us to model aspects such as default in financial applications. We thus extend the result in~\cite{JonteSFDE} on the one hand by considering a more general driving noise but also by considering both the infinite and the random horizon settings. Our treatment of the random horizon problem also motivates the exploration of partial information as optimal controls may be fundamentally different in the partial information setting.

The main contributions of the present work are twofold. First we show that the problem of maximizing $J$ has a solution under certain assumptions on $\varphi$ and $c$, summarized in the definition of an \emph{admissible reward pair}, by finding an optimal control in terms of a family of interconnected value processes. We refer to this family of processes as a \emph{verification family}. Furthermore, we give a set of conditions under which the reward pair defined by \eqref{ekv:phiSFDE}-\eqref{ekv:cSFDE} is admissible.

The remainder of the article is organized as follows. In the next section we state the problem, set the notation used throughout the article and detail the set of assumptions that are made. In particular we introduce the notion of an admissible reward pair. Furthermore, we recall some results, such as the Snell envelope, that will be useful when showing existence of optimal controls. Then, in Section~\ref{sec:VERthm} a verification theorem is derived. This verification theorem is an extension of the verification theorem for the multi-modes optimal switching problem with memory developed in~\cite{SwitchElephant} and presumes the existence of a verification family. In Section~\ref{sec:exist} we show that, under the assumptions made, there exists a verification family whenever $(\varphi,c)$ is an admissible reward pair, thus proving existence of an optimal control for the impulse control problem with the cost functional $J$ defined in \eqref{ekv:objfun}. Then, in Section~\ref{sec:SFDEs} we show that a type of impulse control problems for controlled SFDEs satisfies the conditions on $\varphi$ and $c$ as prescribed in the definition of an admissible reward pair, both in the infinite and random horizon settings.

\section{Preliminaries}
We let $(\Omega,\mcG,\Prob)$ be a complete probability space, and $\bbF:=(\mcF_t)_{t\geq 0}$ a filtration of sub-$\sigma$-fields of $\mcG$ satisfying the usual conditions in addition to being quasi-left continuous and define $\mcF:=\mcF_\infty$.

\begin{rem}
Recall here the concept of quasi-left continuity: A \cadlag process $(X_t:t\geq 0)$ is quasi-left continuous if for each predictable stopping time $\theta$ and every announcing sequence of stopping times $\theta_k\nearrow\theta$ we have $X_{\theta -}:=\lim\limits_{k\to\infty}X_{\theta_k} = X_\theta$, $\Prob$-a.s. Similarly, $X$ is quasi-left upper semi-continuous if $X_{\theta -}\leq  X_\theta$, $\Prob$-a.s. A filtration is quasi-left continuous if $\mcF_{\theta}=\mcF_{\theta-}$ for every predictable stopping time $\theta$.
\end{rem}

\noindent Throughout, we will use the following notation:
\begin{itemize}
  \item $\mcP_{\bbF}$ is the $\sigma$-algebra of $\bbF$-progressively measurable subsets of $\R_+\times \Omega$.
  \item For $p\geq 1$, we let $\mcS^{p}$ be the set of all $\R$-valued, $\mcP_{\bbF}$-measurable, \cadlag processes $(Z_t: t\geq 0)$ such that $\E\left[\sup_{t\in\R_+} |Z_t|^p\right]<\infty$ and let $\Sqlc^{p}$ be the subset of processes that are quasi-left upper semi-continuous. 
  \item We let $\mcT$ be the set of all $\bbF$-stopping times and for each $\eta\in\mcT$ we let $\mcT_\eta$ be the corresponding subsets of stopping times $\tau$ such that $\tau\geq \eta$, $\Prob$-a.s. Furthermore, we let $\mcT^f$ (resp. $\mcT^f_\eta$) be the subset of $\mcT$ (resp. $\mcT_\eta$) with all stopping times $\tau$ for which $\Prob[\tau<\infty]=1$.
  \item For each $\tau\in\mcT$, we let $\mcA(\tau)$ be the set of all $\mcF_\tau$-measurable random variables taking values in $U$.
  \item We let $\mcU$ be the set of all $u=(\tau_1,\ldots,\tau_N;\beta_1,\ldots,\beta_N)$, where $(\tau_j)_{j=1}^N$ is a non-decreasing sequence of $\bbF$-stopping times and $\beta_j\in\mcA(\tau_j)$.
  \item We let $\mcU^f$ denote the subset of $u\in\mcU$ for which $N(t):=\sup\{j:\tau_j\leq t\}$ is $\Prob$-a.s.~finite on compacts (\ie $\mcU^f:=\{u\in\mcU:\: \Prob\left[\{\omega\in\Omega : N(t)>k, \:\forall k>0\}\right]=0, \:\forall t\in\R_+\}$) and for all $k\geq 0$ we let $\mcU^k:=\{u\in\mcU:\:N\leq k,\,\Prob{\rm - a.s.}\}$.
  \item For a random interval $A$ (\ie a set of the type $[\eta_1,\eta_2]$, $[\eta_1,\eta_2)$, $(\eta_1,\eta_2]$ or $(\eta_1,\eta_2)$ for some $\eta_1,\eta_2\in\mcT$) we let $\mcU_{A}$ (and $\mcU_{A}^f$ resp.~$\mcU_{A}^k$) be the subset of $\mcU$ (and $\mcU^f$ resp.~$\mcU^k$) with $\tau_j\in A$, $\Prob$-a.s.~for $j=1,\ldots,N$. When the interval is $A=[\eta,\infty]$ for some $\eta\in\mcT$ we use the shorthand $\mcU_{\eta}$ (and $\mcU_{\eta}^f$ resp.~$\mcU_{\eta}^k$).
  \item We let $\mcD^f$ be the subset of $\mcD$ with all finite sequences and for $k\geq 0$ we let $\mcD^k:=\cup_{i=0}^k D^i$.
  \item Throughout, we let $\vecv=(\vect,\vecb)$, with $\vect:=(t_1,\ldots,t_n)$ and $\vecb:=(b_1,\ldots,b_n)$, where $n$ is possibly infinite, denote a generic element of $\mcD$.
  \item For $\vecv=(\vect,\vecb)\in\mcD^f$ and $\vecv'=(\vect',\vecb')\in\mcD$ we introduce the composition, denoted by $\circ$, defined as $\vecv\circ\vecv':=(t_1,\ldots,t_n,t'_1\vee t_n,\ldots,t'_{n'}\vee t_n;b_1,\ldots,b_n,b'_1,\ldots,b'_{n'})$. Furthermore, for $\vecv\in\mcD$, we define the truncation to $k\geq 0$ interventions as $[\vecv]_{k}:=(t_1,\ldots,t_{k\wedge n};b_1,\ldots,b_{k\wedge n})$
  \item For $l\geq 0$, we let $\Pi_l:=\{0,1/2^{l},2/2^{l},\ldots\}$ and set $\bold \Pi:=\cup_{l=1}^\infty \Pi_l$.
\end{itemize}
Furthermore, we define the following set:
\begin{defn}\label{def:mcH}
We let $\mcH'_\bbF$ be the set of all $\mcP_\bbF\otimes\mcB(U)$-measurable maps\footnote{Throughout, we generally suppress dependence on $\omega$ and refer to $h\in\mcH_\bbF$ as a map $(t,b)\to h(t,b)$.} $h:\Omega\times\R_+\times U\to\R$ such that the collection $\{h(\tau,\beta): \tau\in\mcT^f,\,\beta\in\mcA(\tau)\}$ is uniformly integrable and (outside of a $\Prob$-null set) we have for all $(t,b)\in \R_+\times U$:
\begin{enumerate}[i)]
  \item\label{mcH:lim} The limit $\lim_{(t',b')\to(t,b)}h(t',b')$ exists,
  \item\label{mcH:cont} $\lim_{t'\searrow t}\sup_{b'\in U}|h(t',b')-h(t,b')|=0$, and
  \item\label{mcH:usc} $\lim_{b'\to b}h(t,b')\leq h(t,b)$.
\end{enumerate}
Furthermore, we let $\mcH_\bbF$ be the set of all $h\in\mcH'_\bbF$ such that for any predictable stopping time $\theta\in\mcT$ and any announcing sequence $\theta_j\nearrow \theta$ with $\theta_j\in\mcT^f$ we have $\limsup_{j\to\infty}\sup_{b\in U}\{h(\theta_j,b)-h(\theta,b)\}\leq 0$, $\Prob$-a.s.
\end{defn}


\subsection{Problem formulation}
In the above notation, our problem is characterized by the complete probability space $(\Omega,\mcG,\Prob)$ and the following objects:
\begin{itemize}
  \item The filtration $\bbF$.
  \item A $\mcG\otimes \mcB(\mcD)$-measurable map $\varphi:\mcD \to\R$.
  \item A $\mcG\otimes \mcB(\mcD)$-measurable map $c:\mcD \to\R_+$.
\end{itemize}

To obtain existence of optimal controls we need to make some assumptions on the involved objects. The assumptions that we will use are summarized in the definition of what we refer to as an admissible reward pair:
\begin{defn}\label{ass:onPSI}
We call the pair $(\varphi,c)$ an \emph{admissible reward pair} if:
\begin{enumerate}[(i)]
  \item\label{ass:psiBND} The terminal reward $\varphi$ and the intervention cost $c$ satisfy the following bounds:
  \begin{enumerate}[a)]
    \item $\sup_{u\in\mcU}\E[ |\varphi(u)|^2]<\infty$.
    \item $\sup_{u\in\mcU}\E[ |c(u)|^2]<\infty$ and $c(\vecv)\geq \delta(t_n)$ for all $\vecv\in\mcD^f$, $\Prob$-a.s., where $\delta: \R_ +\to\R_ +$ is a deterministic, continuous, non-increasing and positive function, \ie $\delta(s)\geq\delta(t)>0$, whenever $0\leq s\leq t<\infty$.
  \end{enumerate}
  \item\label{ass:psiREG} For every $v\in\mcU^f$ and every $k\geq 0$ and $T>0$, there are maps $\chi^T\in\mcH'_\bbF$, $\chi\in\mcH_\bbF$ and $(t,b)\mapsto -c^v_{t,b}\in\mcH_\bbF$ such that for all $\tau\in \mcT$ and $b\in U$ we have
  \begin{align*}
    \chi(\tau,b)&=\esssup_{u\in\mcU^k_\tau}\big[\varphi(v\circ(\tau,b)\circ u)-\sum_{j=1}^N c(v\circ(\tau,b)\circ[u]_j)\big|\mcF_\tau\big],
    \\
    \chi^T(\tau,b)&=\esssup_{u\in\mcU^k_{[\tau,T)}}\big[\varphi(v\circ(\tau,b)\circ u)-\sum_{j=1}^N c(v\circ(\tau,b)\circ[u]_j)\big|\mcF_\tau\big],
  \end{align*}
  and
  \begin{align*}
    c^v_{\tau,b}=\big[c(v\circ(\tau,b))\big|\mcF_\tau\big],
  \end{align*}
  $\Prob$-a.s.~(with an exception set that is independent of $b$).
    \item\label{ass:@end} We have $\sup_{u\in\mcU,v\in\mcU_T}\E\big[|\varphi(u\circ v)-\varphi(u)|^2\big]\to 0$ as $T\to\infty$
\end{enumerate}
\end{defn}


The conditions in the above definition are mainly standard assumptions for infinite horizon stochastic impulse control problems translated to our setting. Condition (\ref{ass:psiBND}.a) together with positivity of the intervention cost, $c$, in (\ref{ass:psiBND}.b) implies that the expected maximal reward is finite. Condition~(\ref{ass:@end}) implies that the future has diminishing impact on the total reward and can be seen as a generalization of the deterministic discounting applied in \eqref{ekv:std-imp-cont}. We show below that the fact that the intervention costs are bounded from below by a positive function together with (\ref{ass:psiBND}.a) implies that, with probability one, the optimal control (whenever it exists) can only make a finite number of interventions within any compact time interval.

\begin{rem}
Note that we may hide part of the intervention cost within the function $\varphi$ which implies that, similar to the setting in~\cite{MartyrSigned}, we can handle problems with negative intervention costs as long as a type of martingale condition is satisfied.
\end{rem}

We consider the following problem:\\

\noindent\textbf{Problem 1.} Find $u^*\in\mcU$, such that
\begin{equation}\label{ekv:OPTprob}
J(u^*)=\sup_{u\in\mcU} J(u),
\end{equation}
when $(\varphi,c)$ is an admissible reward pair.\qed\\

Throughout Section \ref{sec:VERthm}-\ref{sec:exist} we will thus assume that $(\varphi,c)$ is an admissible reward pair, before we in Section~\ref{sec:SFDEs} give a set of conditions under which we are able to show that a particular $(\varphi,c)$ of the form \eqref{ekv:phiSFDE}-\eqref{ekv:cSFDE} is an admissible reward pair.

As a step in solving Problem 1 we need the following proposition which is a standard result for impulse control problems.
\begin{prop}\label{prop:finSTRAT} Suppose that there is a $u^*\in\mcU$ such that $J(u^*)\geq J(u)$ for all $u\in\mcU^f$. Then $u^*$ is an optimal control for Problem 1, \ie $J(u^*)\geq J(u)$ for all $u\in\mcU$.
\end{prop}

\noindent\emph{Proof.} Pick $\hat u:=(\hat\tau_1,\ldots,\hat\tau_{\hat N};\hat\beta_1,\ldots,\hat\beta_{\hat N})\in \mcU\setminus \mcU^f$. Then there is a $t\in\R_+$ such that $\Prob[B]>0$ with $B:=\{\omega \in\Omega: \hat N(t)>k, \:\forall k>0\}$. Furthermore, by positivity of the intervention costs\footnote{Throughout $C$ will denote a generic positive constant that may change value from line to line.}
\begin{align*}
J(\hat u)&\leq \sup_{u\in\mcU} \E\big[ |\varphi(u)|\big]-k\delta(t)\Prob[B]\leq C-k\delta(t)\Prob[B],
\end{align*}
for all $k\geq 0$, by Definition~\ref{ass:onPSI}.(\ref{ass:psiBND}.a). However, again by Definition~\ref{ass:onPSI}.(\ref{ass:psiBND}.a) we have $J(\emptyset)\geq -C$. Hence, $\hat u$ is dominated by the strategy of doing nothing and the assertion follows.\qed


\subsection{The Snell envelope}
In this section we gather the main results concerning the Snell envelope that will be useful later on. Recall that a progressively measurable process $X$ is of class [D] if the set of random variables $\{X_\tau:\tau\in\mcT^f\}$ is uniformly integrable.

\begin{thm}[The Snell envelope]\label{thm:Snell}
Let $X=(X_t)_{t\geq 0}$ be an $\bbF$-adapted, $\R$-valued, \cadlag process of class [D]. Then there exists a unique (up to indistinguishability), $\R$-valued \cadlag process $Z=(Z_t)_{t\geq 0}$ called the Snell envelope of $X$, such that $Z$ is the smallest supermartingale that dominates $X$. Moreover, the following holds (with $\Delta X_t:=X_{t}-X_{t-}$):
\begin{enumerate}[(i)]
  \item\label{Snell:sup} For any stopping time $\eta$,
    \begin{equation}\label{ekv:SnellZ}
      Z_{\eta}=\esssup_{\tau\in \mcT_{\eta}}\E\left[X_\tau\big|\mcF_\eta\right].
    \end{equation}
  \item\label{Snell:DoobMeyer} The Doob-Meyer decomposition of the supermartingale $Z$ implies the existence of a triple $(M,K^c,K^d)$ where $(M_t:t\geq 0)$ is a uniformly integrable right-continuous martingale, $(K^c_t:t\geq 0)$ is a non-decreasing, predictable, continuous process with $K^c_0=0$ and $(K^d_t:t\geq 0)$ is non-decreasing purely discontinuous predictable with $K^d_0=0$, such that
      \begin{equation}\label{ekv:DoobMeyerDec}
        Z_t=M_t-K^c_t-K^d_t.
      \end{equation}
      Furthermore, $\{\Delta K_t^d>0\}\subset \{\Delta X_t<0\}\cap\{Z_{t-}=X_{t-}\}$ for all $t\geq 0$.
  \item\label{Snell:att} Let $\eta\in\mcT$ be given and assume that for any predictable $\theta\in\mcT_\eta$ and any increasing sequence $\{\theta_j\}_{j\geq 0}$ with $\theta_j\in\mcT^f_\eta$ and $\lim_{j\to\infty}\theta_j=\theta$, $\Prob$-a.s, we have 
$\limsup_{j\to\infty}X_{\theta_j}\leq X_{\theta}$, $\Prob$-a.s. Then, the stopping time $\tau^*_{\eta}$ defined by $\tau^*_{\eta}:=\inf\{s\geq\eta:Z_s=X_s\}$ (with the convention that $\inf\emptyset=\infty$) is optimal after $\eta$, \ie
    \begin{equation*}
      Z_{\eta}=\E\left[X_{\tau^*_\eta}\big|\mcF_\eta\right].
    \end{equation*}
    Furthermore, in this setting the Snell envelope, $Z$, is quasi-left continuous, \ie $K^d\equiv 0$.
  \item\label{Snell:lim} Let $X^k$ be a sequence of \cadlag processes converging increasingly and pointwisely to the \cadlag process $X$ and let $Z^k$ be the Snell envelope of $X^k$. Then the sequence $Z^k$ converges increasingly and pointwisely to a process $Z$ and $Z$ is the Snell envelope of $X$.
\end{enumerate}
\end{thm}
In the above theorem, (\ref{Snell:sup})-(\ref{Snell:att}) are standard results and proofs can be found in, for example, \cite{ElKarouiLN,HamRefBSDE}. A finite horizon version of statement (\ref{Snell:lim}), which extends trivially to infinite horizon, was proved in \cite{BollanMSwitch1}.\\

The Snell envelope will be the main tool in showing that Problem 1 has a solution.


\subsection{The section and projection theorems}
In this section we recall two fundamental results from the general theory of stochastic processes, namely the measurable selection and the optional projection theorems.

For any space $E$ we define the projection of a set $A\subset \Omega\times E$ onto $\Omega$ as $\pi_\Omega(A):=\{\omega\in\Omega:\exists x\in E,\,(\omega,x)\in A\}$.

\begin{thm}[Measurable projection]\label{thm:MeasProj}
Let $(\Omega,\mcF,\Prob)$ be a complete probability space and $E$ a Polish space. For every $A\in \mcF\otimes \mcB(E)$ the set $\pi_\Omega(A)$ is $\mcF$-measurable.
\end{thm}
A proof can be found in \eg \cite{DelMeyer1} Chapter III. In particular we need the following corollary result:
\begin{cor}\label{cor:MPsup}
Consider a complete probability space $(\Omega,\mcF,\Prob)$ and let $h(\omega,x)$ be a real valued, measurable function defined on the product space $(\Omega\times \R^m,\mcF\otimes\mcB(\R^m))$. Then for all $A\in \mcF\otimes\mcB(\R^m)$, the function
\begin{align*}
g(\omega):=\sup_{x\in\R^m}\{h(\omega,x):(\omega,x)\in A\}
\end{align*}
(with the convention $\sup\emptyset=-\infty$) is $\mcF$-measurable.
\end{cor}
\noindent\emph{Proof.} For each $K\in \R$ we have $\{g(\omega)> K\}=\pi_\Omega(A\cap h^{-1}((K,\infty]))$. Now, since $h$ is measurable, the set $A\cap h^{-1}((K,\infty])$ is in $\mcF\otimes\mcB(\R^m)$ and the result follows by the measurable projection theorem.\qed\\

\begin{thm}[Measurable selection]\label{thm:MeasSel}
Let $(\Omega,\mcF,\Prob)$ be a complete probability space and let $(E,\mcE)$ be a Borel space with $\mcE:=\mcB(E)$. For every $A\in \mcF\otimes \mcB(E)$ there is a $\mcF$-measurable r.v.~$\beta$ taking values in $\bar E:=E\cup\{\partial\}$ (with $\partial$ a cemetery point) such that
\begin{align*}
\{(\omega,\beta(\omega))\in\Omega\times E\}\subset A \quad {\rm and}\quad \{\omega\in\Omega:\beta(\omega)\in E\}=\pi_\Omega(A).
\end{align*}
\end{thm}
This is a standard result and proofs can be found in \eg Chapter 7 in \cite{BertsekasShreve} or \cite{ElKarouiTan}. In particular we need the following well known corollary result:
\begin{cor}\label{cor:MSattained}
Consider a complete probability space $(\Omega,\mcF,\Prob)$ and let $h(\omega,x)$ be a measurable function defined on the product space $(\Omega\times \R^m,\mcF\otimes\mcB(\R^m))$, such that for $\Prob$-almost every $\omega$ the map $x\mapsto h(\omega,x)$ is upper semi-continuous. Then, with $U$ a compact subset of $\R^m$, there exists a $\mcF$-measurable r.v.~$\beta$ such that
\begin{align*}
h(\omega,\beta(\omega))=\sup_{x\in\R^n}\{h(\omega,x):(\omega,x)\in \Omega\times U\},
\end{align*}
$\Prob$-a.s.
\end{cor}
\noindent\emph{Proof.} Since $A:=\Omega\times U\in \mcF\otimes \mcB(E)$ (where now $E=\R^m$) the function $g(\omega)=\sup_{x\in E}\{h(\omega,x):(\omega,x)\in A\}$ is $\mcF$-measurable. Furthermore, as $h$ is $\mcF\otimes\mcB(E)$-measurable, the set $B:=\{(\omega,x)\in\Omega\times U:h(\omega,x)=g(\omega)\}$ is in $\mcF\otimes\mcB(E)$. Now, by Theorem~\ref{thm:MeasSel} there is a $\mcF$-measurable $\bar E$-valued r.v.~$\beta$ such that $\{(\omega,\beta(\omega))\in\Omega\times E\}\subset B$ and $\{\omega\in\Omega:\beta(\omega)\in E\}=\pi_\Omega(B)$. As $U$ is compact and $b\mapsto h(\omega,b)$ is u.s.c.~on $\Omega\setminus\mcN$ with $\Prob(\mcN)=0$, we have that $B^\omega:=\{b\in U: (\omega,b)\in B\}=\{b\in U: h(\omega,b)=g(\omega)\}\neq\emptyset$ for all $\omega\in \Omega\setminus\mcN$ and, hence, $\Prob(\pi_\Omega(B))=1$.\qed\\

The last result that we need is the optional projection theorem.
\begin{thm}[Optional projection]\label{thm:OP}
Assume that $(X_t:t\geq 0)$ is a measurable process (not necessarily adapted to the filtration $\bbF$) with $\E\big[|X_\tau|\big]<\infty$ for all stopping times $\tau\in\mcT$, then there exists a unique optional process $(^o\!X_t:t\geq 0)$ such that
\begin{align*}
\ett_{[\tau<\infty]}{}^o\!X_\tau=\E\big[\ett_{[\tau<\infty]}X_\tau\big|\mcF_\tau],
\end{align*}
for all stopping times $\tau\in\mcT$. If, furthermore, $X$ is \cadlag then ${}^o\!X$ is also \cadlagp
\end{thm}
A proof of Theorem~\ref{thm:OP} can be found in Chapter VI of~\cite{DelMeyer2}.

\subsection{Relevant properties of $\mcH_\bbF$}
We note the following properties:
\begin{lem}\label{lem:mcH-prop}
\begin{enumerate}[a)]
  \item If $h,h'\in\mcH_\bbF$ ($h,h'\in\mcH'_\bbF$), then $h+h'$ and $h\vee h'$ are also in $\mcH_\bbF$ ($\mcH'_\bbF$).
  \item If $h\in\mcH'_\bbF$ then there is a $\mcP_\bbF$-measurable \cadlag process, $h^*$, of class [D], such that $h^*_\tau=\sup_{b\in U}h(\tau,b)$, $\Prob$-a.s.~for each $\tau\in\mcT^f$. If $h\in\mcH_\bbF$ then $h^*$ is quasi-left upper semi-continuous.
  \item If $h,h'\in\mcH'_\bbF$, then $(\sup_{b\in U}|h'(t,b)-h(t,b)|:t\geq 0)$ is $\mcP_\bbF$-measurable and \cadlagp
  \item If $(h^k)_{k\geq 0}$ is a  sequence in $\mcH_\bbF$ that converges uniformly to some $h$ (outside of a $\Prob$-null set) then $h\in\mcH_\bbF$.
\end{enumerate}
\end{lem}

\noindent\emph{Proof.} The first property is trivial.

Moving on to property \emph{b)}, we let $s^l_j:=j2^{-l}$ and note that for $j\geq 0$, there is a $\beta^l_{j+1}\in\mcA(s^l_j)$ such that
\begin{align*}
\sup_{b\in U}\E\big[h(s^l_{j+1},b)|\mcF_{s^l_j}\big]=\E\big[h(s^l_{j+1},\beta^l_{j+1})|\mcF_{s^l_j}\big]
\end{align*}
$\Prob$-a.s., by Corollary~\ref{cor:MSattained}. By Theorem~\ref{thm:OP} we can define the sequence of \cadlag processes $(\hat h^{l})_{l\geq 0}$ as
\begin{align*}
\hat h^{l}_t:=\sum_{j=0}^{\infty}\ett_{[s^l_j,s^l_{j+1})}(t)\E\big[h(t,\beta^l_{j+1})|\mcF_{t}\big].
\end{align*}

Now, we let $\underline {h}^{0}:=\hat h^{0}$ and then recursively define $\underline {h}^{l}:=\underline {h}^{l-1}\vee \hat h^{l}$ for $l\geq 1$. Then, $(\underline {h}^{l})_{l\geq 0}$ is a non-decreasing sequence of \cadlag processes and
\begin{align*}
\sup_{b\in U}h(t,b)\geq \underline {h}^{l}_t,
\end{align*}
$\Prob$-a.s., for all $t\in[0,\infty)$ and all $l\geq 0$. Furthermore, for $t\geq 0$ we let $\iota_l:=\max\{j:j2^{-l}\leq t\}$ and get that
\begin{align*}
\sup_{b\in U}h(t,b) - \underline {h}^{l}_t&=\sup_{b\in U}h(t,b)- \E\big[h(s^l_{\iota_l+1},\beta^l_{\iota_l+1})\big|\mcF_{s^l_{\iota_l}}\big] +\E\big[h(s^l_{\iota_l+1},\beta^l_{\iota_l+1})\big|\mcF_{s^l_{\iota_l}}\big] - \E\big[h(t,\beta^l_{\iota_l+1})\big|\mcF_{t}\big]
\\
&\leq \sup_{b\in U}\{h(t,b)- \E\big[h(s^l_{\iota_l+1},b)\big|\mcF_{s^l_{\iota_l}}\big]\}+\E\big[h(s^l_{\iota_l+1},\beta^l_{\iota_l+1})\big|\mcF_{s^l_{\iota_l}}\big] - \E\big[h(t,\beta^l_{\iota_l+1})\big|\mcF_{t}\big].
\end{align*}
Since $s^l_{\iota_l}\nearrow t$ we have, by quasi-left continuity of the filtration and uniform integrability that
\begin{align*}
\lim_{l\to\infty}\E\big[h(s^l_{\iota_l+1},\beta^l_{\iota_l+1})\big|\mcF_{s^l_{\iota_l}}\big] =\lim_{l\to\infty}\E\big[h(s^l_{\iota_l+1},\beta^l_{\iota_l+1})\big|\mcF_{t}\big].
\end{align*}
Now, as $s^l_{\iota_l+1}\searrow t$, it follows by Definition~\ref{def:mcH}.(\ref{mcH:cont}) that $\underline {h}^{l}_t\nearrow \sup_{b\in U}h(t,b)$, $\Prob$-a.s., as $l\to\infty$. We note that $\underline {h}^l$ is an increasing, uniformly bounded, sequence of $\mcP_\bbF$-measurable \cadlag processes. The sequence, thus, converges to a $\mcP_\bbF$-measurable process, $h^*$. It remains to show that $h^*$ is \cadlagSTOP, quasi-left upper semi-continuous and that it agrees with $\sup_{b\in U}h(t,b)$ on stopping times.

To show that the limit is \cadlag we let
\begin{align*}
\bar h^{l}_t:=\sum_{j=0}^{\infty}\ett_{[s^l_j,s^l_{j+1})}(t)\sup_{(r,b)\in [t,s^l_{j+1})\times U}h(r,b).
\end{align*}
We note that $\bar h^l$ has left and right limits and that, furthermore, $\lim_{t'\searrow t}\bar h^l_{t'}\leq \bar h^l_{t}$. Now, if $\lim_{t'\searrow t}\bar h^l_{t'}< \bar h^l_{t}$ then
\begin{align*}
\bar h^l_{t}=\sup_{b\in U}h(t,b)=h(t,\beta_t)
\end{align*}
for some $\mcF_t$-measurable $\beta_t$, but then we would have
\begin{align*}
h(t,\beta_t)>\lim_{t'\searrow t}\bar h^l_{t'}\geq \lim_{t'\searrow t}h(t',\beta_t)
\end{align*}
contradicting the fact that $(h(s,\beta_t):s\geq t)$ is right continuous. We conclude that $(\bar h^l)_{l\geq 0}$ is a non-increasing sequence with $\bar h^l$ a $\{\mcF_{t+2^{-l}}\}_{t\geq 2^{-l}}$-adapted \cadlag process and
\begin{align*}
\sup_{b\in U}h(t,b)\leq \bar{h}^l_t.
\end{align*}
For $t\geq 0$ we know that there is a non-increasing sequence $(\tilde \tau_l)_{l\geq 0}$, with $\tilde\tau_l$ a $\mcF_{t+2^{-l}}$-measurable r.v.~taking values in $[t,t+2^{-l})$ such that
\begin{align*}
\lim_{l\to\infty}\bar{h}^l_t=\lim_{l\to\infty} \sup_{b\in U}h_{\tilde \tau_l}^{\tilde \tau_l,b}.
\end{align*}
Now, as $h\in\mcH_\bbF$ we have $\lim_{t'\searrow t}\sup_{b\in U}|h(t',b)-h(t,b)|= 0$. It follows that
\begin{align*}
\lim_{l\to\infty}\bar{h}^l_t-\sup_{b\in U}h(t,b)=\lim_{l\to\infty}\sup_{b\in U}h(\tilde \tau_l,b)-\sup_{b\in U}h(t,b)\leq \lim_{l\to\infty}\sup_{b\in U}\{h(\tilde \tau_l,b)-h(t,b)\}= 0,
\end{align*}
$\Prob$-a.s., and in particular we find that $\bar{h}^l_t\searrow h^*_t$, $\Prob$-a.s.~as $l\to\infty$. This gives that for any sequence $t_j\searrow t$
\begin{align*}
\liminf_{t_j\searrow t}h^*_{t_j}&\geq \liminf_{t_j\searrow t}\underline {h}^l_{t_j}=\underline {h}^l_{t}
\end{align*}
and
\begin{align*}
\limsup_{t_j\searrow t}h^*_{t_j}&\leq \limsup_{t_j\searrow t}\bar h^l_{t_j}=\bar h^l_{t}.
\end{align*}
Letting $l$ tend to infinity we find that
\begin{align*}
\liminf_{t_j\searrow t}h^*_{t_j}=\limsup_{t_j\searrow t}\sup_{b\in U}h^*_{t_j}=h^*_{t}.
\end{align*}
Similarly we get existence of left-limits and by the uniform integrability property imposed on members of $\mcH_\bbF$ we conclude that $h^*$ is a \cadlagSTOP, $\mcP_\bbF$-measurable process of class [D].\\

Let $(\tau_l)_{l\geq 0}$ be a non-increasing sequence of stopping times in $\mcT^{\bold \Pi}$ (the subset of $\mcT$ with all stopping times taking values in the countable set $\bold \Pi$) such that $\tau_l\searrow\tau$. We may, for example, set $\tau_l:=\inf\{s\in\Pi_l:s\geq \tau\}$. Since $\bold \Pi$ is countable we have
\begin{align}\label{ekv:Y0stop2}
h^*_{\tau_l}=\sup_{b\in U} h(\tau_l,b),
\end{align}
$\Prob$-a.s. Now, by right-continuity we get that $\lim_{l\to\infty}h^*_{\tau_l}=h^*_{\tau}$. Letting $(\beta_l)_{l\geq 0}$ be a sequence of maximizers for the right-hand side of \eqref{ekv:Y0stop2} at times $(\tau_l)_{l\geq 0}$ we get
\begin{align*}
\limsup_{l\to\infty}\sup_{b\in U}h(\tau_l,b)&=\limsup_{l\to\infty}h(\tau_l,\beta_l),
\end{align*}
$\Prob$-a.s. Moreover, for each $\omega\in\Omega$ there is a subsequence $(\iota_j(\omega))_{j\geq 1}$ such that
\begin{align*}
\limsup_{l\to\infty}h(\tau_l,\beta_l) = \lim_{j\to\infty}h(\tau_{\iota_j},\beta_{\iota_j}).
\end{align*}
Since $U$ is compact, there is a subsequence $(\iota'_j(\omega))_{j\geq 1}\subset (\iota_j(\omega))_{j\geq 1}$ such that $(\beta_{\iota'_j}(\omega))_{j\geq 0}$ converges to some $\tilde\beta(\omega)\in U$ and so we have
\begin{align*}
\limsup_{l\to\infty}\sup_{b\in U}h(\tau_l,b)&= \lim_{j\to\infty}(h(\tau_{\iota'_j},\beta_{\iota'_j})-h(\tau,\beta_{\iota'_j}))
 + \lim_{j\to\infty}(h(\tau,\beta_{\iota'_j})-h(\tau,\tilde\beta)) + h(\tau,\tilde\beta).
\end{align*}
Now,
\begin{align*}
\lim_{j\to\infty}(h(\tau_{\iota'_j},\beta_{\iota'_j})-h(\tau,\beta_{\iota'_j}))=0,
\end{align*}
$\Prob$-a.s., by Definition~\ref{def:mcH}.(\ref{mcH:cont}) and
\begin{align*}
\lim_{j\to\infty}(h(\tau,\beta_{\iota'_j})-h(\tau,\tilde\beta))\leq 0,
\end{align*}
$\Prob$-a.s., by the upper semi-continuity declared in Definition~\ref{def:mcH}.(\ref{mcH:usc}). Further, as $h(\tau,\tilde\beta)\leq \sup_{b\in U}h(\tau,b)$ we conclude that
\begin{align*}
\limsup_{l\to\infty}\sup_{b\in U}h(\tau_l,b)\leq\sup_{b\in U}h(\tau,b),
\end{align*}
$\Prob$-a.s. On the other hand, there is a $\beta\in\mcA(\tau)$ such that $\sup_{b\in U}h(\tau,b)=h(\tau,\beta)$, $\Prob$-a.s., and we have
\begin{align*}
\liminf_{l\to\infty}\sup_{b\in U}h(\tau_l,b)&=\liminf_{l\to\infty}h(\tau_l,\beta_l) \geq \lim_{l\to\infty}h(\tau_l,\beta)
=h(\tau,\beta),
\end{align*}
$\Prob$-a.s. This implies that the limit exists with $\lim_{l\to\infty}\sup_{b\in U}h(\tau_l,b)=\sup_{b\in U}h(\tau,b)$, $\Prob$-a.s., establishing that $h^*_\tau=\sup_{b\in U}h(\tau,b)$, $\Prob$-a.s. Finally, quasi-left upper semi-continuity of $h^*$ is immediate from Definition~\ref{def:mcH} and b) follows.

Property c) follows similarly by noting that for any $\epsilon_l>0$ we can chose $\beta^l_{j+1}\in\mcA(s^l_j)$ such that
\begin{align*}
\sup_{b\in U}\E\big[|h'(s^l_{j+1},b)-h(s^l_{j+1},b)|\big|\mcF_{s^l_j}\big] =\E\big[|h'(s^l_{j+1},\beta^l_{j+1})-h(s^l_{j+1},\beta^l_{j+1})|\big|\mcF_{s^l_j}\big]+\epsilon_l
\end{align*}
and we can choose the sequence $(\epsilon_l)_{l\geq 0}$ such that $\epsilon_l\searrow 0$.

Concerning the last property we note that for each $\epsilon>0$ we can, by uniform convergence, choose a $\Prob$-a.s.~finite $k(\omega)\geq 0$ such that $|h(t,b)-h^k(t,b)|\leq\epsilon$ for all $(t,b)\in\R_+\times U$. Then,
\begin{align*}
\limsup_{(t',b')\to (t,b)}h(t',b')-\liminf_{(t',b')\to (t,b)}h(t',b')\leq \limsup_{(t',b')\to (t,b)}h^k(t',b')-\liminf_{(t',b')\to (t,b)}h^k(t',b')+2\epsilon
=2\epsilon
\end{align*}
and property i) in Definition~\ref{def:mcH} follows as $\epsilon>0$ was arbitrary.
The remaining properties follow similarly.\qed\\

\section{A verification theorem\label{sec:VERthm}}
Our approach to finding a solution to Problem 1 is based on deriving an optimal control under the assumption that a specific family of processes exists, and then showing that the family does indeed exist. We will refer to any such family of processes as a \emph{verification family}. Before making precise the concept of a verification family we introduce the notion of consistency:
\begin{defn}
We refer to a family of processes $((X^v_t)_{t\geq 0}:v\in \mcU^f)$ as being \emph{consistent} if for each $u\in\mcU$, the map $h:\Omega\times\R_+\times U\to\R$ given by $h(t,b)=X^{u\circ(t,b)}_t$ is $\mcP_\bbF\otimes\mcB(U)$-measurable and for each $\tau$ and each $\beta\in\mcA(\tau)$ we have $X^{u\circ(\tau,\beta)}_\tau=h(\tau,\beta)$, $\Prob$-a.s.
\end{defn}

We are now ready to state the definition of a verification family:
\begin{defn}\label{def:vFAM}
We define a \emph{verification family} to be a consistent family of \cadlag supermartingales $((Y^{v}_s)_{s\geq 0}: v\in \mcU^f)$ such that for each $v\in\mcU^f$:
\begin{enumerate}[a)]
  \item\label{vfass:recur} The family satisfies the recursion
  \begin{align}
  Y^{v}_s&=\esssup_{\tau \in \mcT_{s}} \E\Big[\ett_{[\tau = \infty]}\varphi(v)+\ett_{[\tau < \infty]}\sup_{b\in U} \{-c^v_{\tau,b}+Y^{v\circ (\tau,b)}_\tau\}\Big| \mcF_s\Big]\label{ekv:Ydef}
  \end{align}
  \item\label{vfass:uniBND} The family is uniformly bounded in the sense that $\sup_{u\in\mcU^f}\E[\sup_{s\in [0,\infty]}|Y^{u}_s|^2]<\infty$.
  \item\label{vfass:in-mcH} The map $(t,b)\mapsto Y^{v\circ(t,b)}_t$ belongs to $\mcH_\bbF$.
  \item\label{vfass:discount} $\sup_{u\in\mcU^f}\E\big[\sup_{s\in[T,\infty]}|Y^{u}_s-\E[\varphi(u)|\mcF_s]|\big]\to 0$, as $T\to \infty$.
\end{enumerate}
\end{defn}

The purpose of the present section is to reduce the solution of Problem 1 to showing existence of a verification family. This is accomplished by the following verification theorem:
\begin{thm}\label{thm:vfc}
Assume that there exists a verification family $((Y^{v}_s)_{s\geq 0}: v\in \mcU^f)$. Then the family is unique (\ie there is at most one verification family, up to indistinguishability of the maps $t\mapsto Y^{v}_t$ and $(t,b)\mapsto Y^{v\circ(t,b)}_t$) and:
\begin{enumerate}[(i)]
  \item Satisfies $Y_0=\sup_{u\in \mcU} J(u)$ (where $Y:=Y^{\emptyset}$).
  \item Defines the optimal control, $u^*=(\tau_1^*,\ldots,\tau_{N^*}^*;\beta_1^*,\ldots,\beta_{N^*}^*)$, for Problem 1, where $(\tau_j^*)_{1\leq j\leq {N^*}}$ is a sequence of $\bbF$-stopping times given by
  \begin{align}
  \tau^*_j:=\inf \Big\{&s \geq \tau^*_{j-1}:\:Y_s^{[u^*]_{j-1}}=\sup_{b\in U} \{-c^{[u^*]_{j-1}}_{s,b}+Y^{[u^*]_{j-1}\circ(s,b)}_s\}\Big\},\label{ekv:taujDEF}
  \end{align}
  using the convention that $\inf\emptyset=\infty$, the sequence $(\beta_j^*)_{1\leq j\leq {N^*}}$ is defined recursively as a measurable selection of
  \begin{equation}\label{ekv:betajDEF}
  \beta^*_j\in\mathop{\arg\max}_{b\in U}\{-c^{[u^*]_{j-1}}_{\tau^*_j,b}+ Y^{[u^*]_{j-1}\circ(\tau^*_j,b)}_{\tau^*_j}\}
  \end{equation}
  and $N^*=\sup\{j:\tau^*_j<\infty\}$, with $\tau_0^*:=0$.
\end{enumerate}
\end{thm}


\noindent\emph{Proof.} The proof is divided into three steps where we first, in Step 1, show that for any $0\leq j\leq N^*$ we have
\begin{align}
Y^{[u^*]_{j}}_{s}&=\E\Big[\ett_{[\tau^*_{j+1} = \infty]}\varphi(u^*)+\ett_{[\tau^*_{j+1} < \infty]}\{-c([u^*]_{j+1})+Y^{[u^*]_{j+1}}_{\tau^*_{j+1}}\}\Big| \mcF_{s}\Big],\label{ekv:esssupATT}
\end{align}
$\Prob$-a.s.~for $s\in[\tau_j^*,\tau^*_{j+1}]$. Following this, in Step 2, we show that $Y_0=J(u^*)$. Then in Step 3 we show that $u^*$ is the optimal control, establishing \emph{(i)} and \emph{(ii)}. A straightforward generalization to arbitrary initial conditions $v\in\mcU^f$ then gives that
\begin{align}
Y^{v}_s=\esssup_{u\in \mcU_{s}} \E\Big[\varphi(v\circ u)-\sum_{j=1}^{N} c(v\circ [u]_{j})\Big|\mcF_s\Big],\label{ekv:Yuni}
\end{align}
by which uniqueness follows.\\

\noindent {\bf Step 1} We start by showing that for each $v\in \mcU^f$ the recursion \eqref{ekv:Ydef} can be written in terms of a $\bbF$-stopping time and that the inner supremum is attained, $\Prob$-a.s. From \eqref{ekv:Ydef} and consistency we note that $Y^{v}$ is the smallest supermartingale that dominates
\begin{align}
R^{v}:=\Big(&\ett_{[s=\infty]}\E\big[\varphi(v)\big|\mcF\big]+\ett_{[s < \infty]}\sup_{b\in U}\{-c^v_{s,b}+Y^{v\circ(s,b)}_s\}\Big| :\: 0\leq s\leq \infty\Big).\label{ekv:dominated}
\end{align}
By property \ref{vfass:in-mcH}) in the definition of a verification family (Definition~\ref{def:vFAM}) and Definition~\ref{ass:onPSI}.(\ref{ass:psiREG}) we have that the map $(t,b)\mapsto -c^v_{t,b}+Y^{v\circ(t,b)}_t$ belongs to $\mcH_\bbF$. Then, by Lemma~\ref{lem:mcH-prop}.(b) and Definition~\ref{ass:onPSI}.(\ref{ass:psiBND}) it follows that $R^{v}$ is a \cadlag process of class [D] that is quasi-left upper semi-continuous on $[0,\infty)$. Furthermore, by property \ref{vfass:discount}) and positivity of the intervention costs we note that $\limsup_{t\to\infty}R^v_t\leq R^v_\infty$, $\Prob$-a.s. By Theorem \ref{thm:Snell}.(\ref{Snell:att}) and consistency it thus follows that for any $\theta\in \mcT^f$, the stopping time $\tau^\theta\in \mcT_{\theta}$ given by
\begin{align*}
\tau^\theta:=\inf\big\{s\geq\theta:Y^{v}_s=\sup_{b\in U}\{-c^v_{s,b}+Y^{v\circ(s,b)}_{s}\}\big\}
\end{align*}
is such that:
\begin{align*}
Y^{v}_\theta=\E\Big[\ett_{[\tau^\theta = \infty]}\varphi(v)+\ett_{[\tau^\theta < \infty]}\sup_{b\in U}\{-c^v_{\tau^\theta,b}+Y^{v\circ(\tau^\theta,b)}_{\tau^\theta}\}\Big| \mcF_\theta\Big].
\end{align*}
Now, since $(t,b)\mapsto -c^v_{t,b}+Y^{v\circ(t,b)}_t\in\mcH_\bbF$ the map $b\mapsto-c^v_{\tau^\theta,b}+Y^{v\circ(\tau^\theta,b)}_{\tau^\theta}$ is $\mcF_{\tau^\theta}\otimes\mcB(U)$-measurable and u.s.c.~on $\{\tau^\theta<\infty\}\setminus\mcN$ for some $\Prob$-null set $\mcN$. Corollary~\ref{cor:MSattained} and consistency then implies that there is a $\beta ^\theta\in\mcA(\tau^\theta)$ such that
\begin{align*}
Y^{v}_\theta=\E\Big[\ett_{[\tau^\theta = \infty]}\varphi(v)+\ett_{[\tau^\theta < \infty]}\{-c(v\circ (\tau^\theta,\beta^\theta))+Y^{v\circ(\tau^\theta,\beta^\theta)}_{\tau^\theta}\}\Big| \mcF_\theta\Big],
\end{align*}
$\Prob$-a.s., and in particular \eqref{ekv:esssupATT} holds.\\

\noindent {\bf Step 2} We now show that $Y_0=J(u^*)$. We start by noting that $Y$ is the Snell envelope of
\begin{align*}
\Big(\ett_{[s=\infty]}\E\big[\varphi(\emptyset)\big|\mcF\big]+\ett_{[s < \infty]}\sup_{b\in U}\{-c_{s,b}+Y^{s,b}_s\} :\: 0\leq s\leq \infty\Big)
\end{align*}
and by Step 1 we thus have
\begin{align*}
Y_0
&=\E\Big[\ett_{[\tau^*_1=\infty]}\varphi(\emptyset)+\ett_{[\tau^*_1 < \infty]}\{-c(\tau^*_1,\beta^*_1)+Y^{\tau^*_1,\beta^*_1}_{\tau^*_1}\}\Big].
\end{align*}

Moving on we pick $j\in\{1,\ldots, N^*\}$ and note that $[u^*]_{j}\in \mcU^f$. But then, by Step 1, we have that
\begin{align*}
Y^{[u^*]_{j}}_{\tau^*_j}
&=\E\Big[\ett_{[\tau^*_{j+1}=\infty]}\varphi(u^*)
+\ett_{[\tau^*_{j+1} < \infty]}\{-c([u^*]_{j+1})+Y^{[u^*]_{j+1}}_{\tau^*_{j+1}}\}\big|\mcF_{\tau^*_j}\Big].
\end{align*}
By induction we get that for each $K\geq 0$ we have,
\begin{align}
Y_0&=\E\Big[\ett_{[N^*\leq K]}\varphi(u^*)-\sum_{j=1}^{N^*\wedge K}c([u^*]_{j})+\ett_{[N^*> K]}\{-c([u^*]_{K+1})+Y^{[u^*]_{K+1}}_{\tau^*_{K+1}}\}\Big].\label{ekv:Y0-trunk}
\end{align}
Now, arguing as in the proof of Proposition~\ref{prop:finSTRAT} and using property \emph{\ref{vfass:uniBND})} we find that $u^*\in\mcU^f$. To show that the right-hand side of \eqref{ekv:Y0-trunk} equals $J(u^*)$ we note that \eqref{ekv:Y0-trunk} can be rewritten as
\begin{align*}
Y_0&=\E\Big[\varphi([u^*]_{K+1})-\sum_{j=1}^{N^*\wedge K+1}c([u^*]_{j})+\ett_{[N^*> K]}\{Y^{[u^*]_{K+1}}_{\tau^*_{K+1}}-\varphi([u^*]_{K+1})\}\Big]
\end{align*}
which gives
\begin{align}\nonumber
|Y_0-J(u^*)|&\leq\E\Big[|\varphi([u^*]_{K+1})-\varphi(u^*)|\Big]
+\E\Big[\sum_{j=K+2}^{N^*}c([u^*]_{j})\Big] \nonumber
\\
&\quad+|\E\Big[\ett_{[N^*> K]}\{Y^{[u^*]_{K+1}}_{\tau^*_{K+1}}-\varphi([u^*]_{K+1})\}\Big]|\label{ekv:Y0-J}
\end{align}
for all $K\geq 0$. For the first term on the right-hand side we note that for any $T>0$ we have
\begin{align*}
&\E\Big[|\varphi([u^*]_{K+1})-\varphi(u^*)|\Big]
\\
&\leq \E\Big[\ett_{[\tau^*_{K+2}<T]}|\varphi([u^*]_{K})-\varphi(u^*)|\Big]+\sup_{u\in\mcU^f,v\in\mcU_T}\E[|\varphi(u\circ v)-\varphi(u)|]
\\
&\leq 2\Prob[\tau^*_{K+2}<T]^{1/2}\sup_{u\in\mcU}\E[|\varphi(u)|^2]^{1/2}+\sup_{u\in\mcU^f,v\in\mcU_T}\E[|\varphi(u\circ v)-\varphi(u)|].
\end{align*}
where we have used H\"older's inequality to arrive at the last inequality. As $(\varphi,c)$ is an admissible reward pair, Definition~\ref{ass:onPSI}.(\ref{ass:@end}) gives that the second term can be made arbitrarily small by choosing $T$ sufficiently large and Definition~\ref{ass:onPSI}.(\ref{ass:psiBND}) implies that the first term tends to zero as $K\to\infty$ for all finite $T$, since $u^*\in\mcU^f$. We, thus, conclude that the first term on the right-hand side in \eqref{ekv:Y0-J} tends to zero as $K\to\infty$.

For the second term we note that letting $K\to\infty$ in \eqref{ekv:Y0-trunk} and using \emph{\ref{vfass:uniBND})} and Definition~\ref{ass:onPSI}.(\ref{ass:psiBND}.a) we find that $\sum_{j=1}^{N^*\wedge l}c([u^*]_{j})$ converges increasingly to a limit in $L^2(\Omega,\Prob)$ as $l\to\infty$. Hence, the second term on the right-hand side of \eqref{ekv:Y0-J} also tends to zero as $K\to\infty$.

Conditioning on $\mcF_{\tau^*_{K+1}}$ in the third term of the right-hand side of \eqref{ekv:Y0-trunk} and noting that $\{N^*> K\}$ is $\mcF_{\tau^*_{K+1}}$-measurable we find that, similar to the above case, we have for any $T\geq 0$ that
\begin{align*}
|\E\Big[\ett_{[N^*> K]}\{Y^{[u^*]_{K+1}}_{\tau^*_{K+1}}-\varphi([u^*]_{K+1})\}\Big]|&\leq\E\Big[|Y^{[u^*]_{K+1}}_{\tau^*_{K+1}}-\E[\varphi([u^*]_{K+1})|\mcF_{\tau^*_{K+1}}]|\Big]
\\
&\leq C\Prob[\tau^*_{K+1}<T]^{1/2}+\sup_{u\in\mcU}\E[\sup_{s\in[T,\infty]}|Y_s^u-\E[\varphi(u)|\mcF_s]|],
\end{align*}
where the second term can be made arbitrarily small by property \emph{\ref{vfass:discount})} and we conclude that $Y_0=J(u^*)$.\\

\noindent {\bf Step 3} It remains to show that the strategy $u^*$ is optimal. 
To do this we pick any other strategy $\hat u:=(\hat\tau_1,\ldots,\hat\tau_{\hat N};\hat\beta_1,\ldots,\hat\beta_{\hat N})\in\mcU^f$. By Step 2 and the definition of $Y_0$ in~\eqref{ekv:Ydef} we have
\begin{align*}
J(u^*)=Y_0 &\geq \E\Big[\ett_{[\hat N = 0]}\varphi(\emptyset) + \ett_{[\hat N > 0]}\sup_{b\in U}\{-c_{\hat \tau_1,b}+Y^{\hat \tau_1;b}_{\hat \tau_1}\}\Big]
\\
&\geq\E\Big[\ett_{[\hat N = 0]}\varphi(\emptyset) + \ett_{[\hat N > 0]}\{-c(\hat \tau_1;\hat\beta_1)+Y^{\hat \tau_1;\hat\beta_1}_{\hat \tau_1}\}\Big]
\end{align*}
but in the same way
\begin{align*}
Y^{\hat \tau_1,\hat \beta_1}_{\hat \tau_1}\geq\E\Big[\ett_{[\hat N = 1]}\varphi(\hat\tau_1,\hat\beta_1) + \ett_{[\hat N > 1]}\{-c(\hat \tau_1,\hat \tau_2;\hat\beta_1,\hat\beta_2)+Y^{\hat \tau_1,\hat \tau_2;\hat\beta_1,\hat\beta_2}_{\hat \tau_2}\}\Big| \mcF_{\hat\tau_1}\Big],
\end{align*}
$\Prob$--a.s. Repeating this procedure $K$ times gives
\begin{align*}
J(u^*)\geq \hat Y_K:=\E\Big[\ett_{[\hat N \leq K]}\varphi(\hat u)-\sum_{j=1}^{\hat N\wedge K}c([\hat u]_{j}) + \ett_{[\hat N > K]}\{-c([\hat u]_{K+1})+Y^{[\hat u]_{K+1}}_{\hat \tau_{K+1}}\}\Big].
\end{align*}
Now, we have
\begin{align*}
J(\hat u)-\hat Y_K&\leq
\E\Big[\ett_{[\hat N>K]}\{\varphi(\hat u)-Y^{[\hat u]_{K+1}}_{\hat \tau_{K+1}}\}\Big]
\\
&= \E\Big[\ett_{[\hat N>K]}\{\varphi([\hat u]_{K+1})-Y^{[\hat u]_{K+1}}_{\hat \tau_{K+1}}\}\Big]
+\E\Big[\ett_{[\hat N>K]}\{\varphi(\hat u)-\varphi([\hat u]_{K+1})\}\Big]
\\
&\leq \E\Big[|\E\big[\varphi([\hat u]_{K+1})\big|\mcF_{\hat \tau_{K+1}}\big]-Y^{[\hat u]_{K+1}}_{\hat \tau_{K+1}}|\Big]
 + \E\Big[|\varphi(\hat u)-\varphi([\hat u]_{K+1})|\Big]
\end{align*}
where the right-hand side tends to zero as $K\to\infty$ by repeating the argument in Step 2, which is possible since $\hat u\in\mcU^f$. We conclude that $J(u^*)\geq J(\hat u)$ for all $\hat u\in\mcU^f$ and it follows by Proposition~\ref{prop:finSTRAT} that $u^*$ is an optimal control for Problem 1.\qed\\


\section{Existence\label{sec:exist}}
Theorem~\ref{thm:vfc} presumes existence of the verification family $((Y^{v}_s)_{s\geq 0}: v\in \mcU^f)$. To obtain a satisfactory solution to Problem 1, we thus need to establish that a verification family exists. This is the topic of the present section. We will follow the standard existence proof which goes by applying a Picard iteration (see \cite{CarmLud,BollanMSwitch1,HamZhang}). We, therefore, let the sequence of consistent families of processes $((Y^{v,k}_s)_{s\geq 0}: v\in \mcU^f)_{k\geq 0}$ satisfy the recursion
\begin{align}
Y^{v,0}_s:=\E\Big[\varphi(v)\Big| \mcF_s\Big],\label{ekv:Y0def}
\end{align}
and
\begin{align}
Y^{v,k}_s:=\esssup_{\tau \in \mcT_{s}} \E\Big[&\ett_{[\tau= \infty]}\varphi(v)
+\ett_{[\tau < \infty]}\sup_{b\in U}\{-c^v_{\tau,b}+Y^{v\circ(\tau,b),k-1}_\tau\}\Big| \mcF_s\Big]\label{ekv:Ykdef}
\end{align}
for $k\geq 1$.\\

We will make use of the following induction hypothesis:

\begin{hyp*}[{\bf VF.k}] There is a sequence of consistent families of \cadlag supermartingales $((Y^{v,k'}_s)_{s\geq 0}: v\in \mcU^f)_{0\leq k'\leq k}$ such that for $k'=0,\ldots,k$ and $v\in \mcU^f$:
\begin{enumerate}[i)]
\item\label{Yk:rec} The relation \eqref{ekv:Y0def} holds for $k'=0$ and \eqref{ekv:Ykdef} holds for $k'>0$.
  \item\label{Yk:in-mcH} The map $(t,b)\mapsto Y^{v\circ(t,b),k'}_t$ belongs to $\mcH_\bbF$.
\end{enumerate}
\end{hyp*}

We note that Hypothesis {\bf VF.k} has similarities with the definition of a verification family. Property (\ref{Yk:rec}) gives a sequential equivalent to \eqref{ekv:Ydef} and property (\ref{Yk:in-mcH}) is equivalent to Definition~\ref{def:vFAM}.(\ref{vfass:in-mcH}). However, we note that the boundedness properties do not appear in Hypothesis {\bf VF.k}. In the following two propositions we show that these are implicit.

\begin{prop}\label{prop:Ykbnd}
Assume that Hypothesis~{\bf VF.k} holds for some $k\geq 0$. Then, the sequence of families of processes $((Y^{v,k'}_s)_{s\geq 0}: v\in \mcU^f)_{0\leq k'\leq k}$ is well defined and uniformly bounded in the sense that there is a $K>0$ (that does not depend on $k$) such that,
\begin{align*}
\sup_{u\in\mcU^f}\E\Big[\sup_{s\in[0,\infty]}|Y^{u,k+1}_s|^2\Big]\leq K.
\end{align*}
Furthermore, for each $v\in\mcD^f$ (whenever it is well defined) the collection $\big\{\sup_{b\in U}\{-c^v_{\tau,b}+Y^{v\circ(\tau,b),k}_\tau\}:\tau\in\mcT^f,\:k\geq 0\big\}$ of random variables is uniformly integrable.
\end{prop}

\noindent\emph{Proof.} We note that under Hypothesis {\bf VF.k}, the sequence of families $((Y^{v,k'}_s)_{s\geq 0}: v\in \mcU^f)_{0\leq k'\leq k+1}$ exists and is uniquely defined up to indistinguishability for each $Y^{v,k'}$ by repeated application of Theorem~\ref{thm:Snell}. By the definition of $Y^{v,k+1}$ and positivity of the intervention costs we have that for any $v\in\mcU^f$,
\begin{align*}
\E\Big[\varphi(v)\big|\mcF_s\Big]\leq Y^{v,k+1}_s&\leq \esssup_{u\in\mcU}\E\Big[\varphi(u)\big|\mcF_s\Big].
\end{align*}
Doob's maximal inequality gives that for any $u\in\mcU$
\begin{align*}
\E\Big[\sup_{s\in[0,\infty]}\E\Big[|\varphi(u)|\big|\mcF_s\Big]^2\Big]\leq C\E\Big[|\varphi(u)|^2\Big].
\end{align*}
Taking the supremum over all $u\in\mcU$ on both sides and using that the right-hand side is uniformly bounded by Definition~\ref{ass:onPSI}.(\ref{ass:psiBND}.a) the first bound follows.\\

Concerning the second claim, note that for each $\tau\in\mcT^f$ and each $\epsilon>0$ repeating the proof of Corollary~\ref{cor:MSattained} with $g^\epsilon:=g-\epsilon$ instead of $g$ we find that there is a $\beta^\epsilon\in\mcA(\tau)$ such that
\begin{align*}
\sup_{b\in U}\{-c^v_{\tau,b}+Y^{v\circ(\tau,b),k}_\tau\} \leq -c^v_{\tau,\beta^\epsilon}+Y^{v\circ(\tau,\beta^\epsilon),k}_\tau+\epsilon
\end{align*}
Now,
\begin{align*}
\E\Big[\sup_{b\in U}|-c^v_{\tau,b}+Y^{v\circ(\tau,b),k}_\tau|^2\Big]
&\leq 2\E\Big[|-c^v_{\tau,\beta^\epsilon}+Y^{v\circ(\tau,\beta^\epsilon),k}_\tau|^2\Big]+2\epsilon^2
\\
&\leq 4\E\Big[\sup_{t\in[0,\infty]}|\E\big[c(v\circ(\tau,\beta^\epsilon))\big|\mcF_t\big]|^2\Big] +4\E\Big[\sup_{t\in[0,\infty]}|Y^{v\circ(\tau,\beta^\epsilon),k}_t|^2\Big]+2\epsilon^2
\end{align*}
where the right-hand side is bounded, uniformly in $(\tau,\beta^\epsilon)$ and $k\geq 0$, by the above in combination with Definition~\ref{ass:onPSI}.(\ref{ass:psiBND}.b) and Doob's maximal inequality.\qed\\

We also have the following diminishing future impact property:
\begin{prop}\label{prop:Ykdiscount}
Assume again that Hypothesis~{\bf VF.k} holds for some $k\geq 0$. Then,
\begin{align*}
\sup_{u\in\mcU}\E\Big[\sup_{t\in[T,\infty]}|Y^{u,k+1}_t-\E\big[\varphi(u)\big|\mcF_t\big]|\Big]\to 0,
\end{align*}
as $T\to\infty$, uniformly in $k$.
\end{prop}

\noindent\emph{Proof.} By the properties of the essential supremum and positivity of the intervention costs we have for every $v\in\mcU^f$,
\begin{align*}
\E\big[\varphi(v)\big|\mcF_t\big]\leq Y^{v,k+1}_t&\leq \esssup_{u\in\mcU_t}\E\big[\varphi(v\circ u)\big|\mcF_t\big].
\end{align*}
This implies that
\begin{align*}
|Y^{v,k+1}_t-\E\big[\varphi(v)\big|\mcF_t\big]|&\leq \esssup_{u\in\mcU_t}\E\big[|\varphi(v\circ u)-\varphi(v)|\big|\mcF_t\big].
\end{align*}
Doob's maximal inequality now gives that for any $u\in\mcU$
\begin{align*}
\E\Big[\sup_{t\in[T,\infty]}|Y^{u,k+1}_t-\E\big[\varphi(u)\big|\mcF_t\big]|^2\Big]&\leq \sup_{u\in\mcU_T}\E\Big[\sup_{s\in [0,\infty]}\E\big[|\varphi(v\circ u)-\varphi(v)|\big|\mcF_t\big]^2\Big]
\\
&\leq C \sup_{u\in\mcU_T}\E\big[|\varphi(v\circ u)-\varphi(v)|^2\big].
\end{align*}
Taking the supremum over all $v\in\mcU^f$ on both sides and using that the right-hand side tends to 0 as $T\to\infty$ by Definition~\ref{ass:onPSI}.(\ref{ass:@end}) the assertion follows by Jensen's inequality.\qed\\

The objective in the remainder of this section is to show that the limit family that we get when letting $k\to\infty$ in $((Y^{v,k}_s)_{s\geq 0}: v\in \mcU^f)$ is a verification family. We start by showing that the induction hypothesis holds:\\

\begin{prop}\label{prop:Yk}
Hypothesis {\bf VF.k} holds for all $k\geq 0$.
\end{prop}

\noindent\emph{Proof.} We start by noting that by Definition~\ref{ass:onPSI}.(\ref{ass:psiREG}) there is a $h\in\mcH_\bbF$ such that for $\tau\in\mcT^f$ we have
\begin{align}\label{ekv:Gam0def}
h(\tau,b)=\E\big[-c(v\circ (\tau,b))+\varphi(v\circ (\tau,b))|\mcF_\tau\big],
\end{align}
for all $b\in U$, $\Prob$-a.s. The statement, thus, holds for $k=0$.

Moving on we assume that {\bf VF.k} holds for some $k\geq 0$. But then, by Propositions~\ref{prop:Ykbnd} and \ref{prop:Ykdiscount} we can applying a reasoning similar to that in the proof of Theorem~\ref{thm:vfc} to find that $Y^{v,k+1}$ is a \cadlag supermartingale with
\begin{align}\label{ekv:Ykalt}
Y^{v,k+1}_t=\esssup_{u \in \mcU^{k+1}_{t}} \E\Big[\varphi(v\circ u)-\sum_{j=1}^{N}c(v\circ [u]_{j})\big|\mcF_t\Big].
\end{align}
By Definition~\ref{ass:onPSI}.(\ref{ass:psiREG}) it follows that there is a consistent family satisfying \eqref{ekv:Ykalt} such that $(t,b)\mapsto Y^{v\circ (t,b),k+1}_t\in\mcH_\bbF$ and we conclude that {\bf VF.k+1} holds as well. By induction this extends to all $k\geq 0$.\qed\\

Having established that Hypothesis {\bf VF.k} holds for all $k\geq 0$ we move on to investigate the limit family that we get when letting $k\to\infty$.
%

\begin{prop}\label{prop:Yklim}
For each $v\in\mcU^f$, the limit $\bar Y^{v}:=\lim_{k\to\infty}Y^{v,k}$, exists as an increasing pointwise limit, $\Prob$-a.s.
\end{prop}

\noindent\emph{Proof.} Since $\mcU^k_t\subset \mcU^{k+1}_t$ we have that, $\Prob$-a.s.,
\begin{align*}
Y^{v,k}_t \leq Y^{v,k+1}_t\leq \esssup_{u\in\mcU_t}\E\Big[|\varphi(v\circ u)|\big|\mcF_t\Big],
\end{align*}
where the right-hand side is bounded $\Prob$-a.s.~by Proposition~\ref{prop:Ykbnd}. Hence, the sequence $(Y^{v,k}_t:t\geq 0)_{k\geq 0}$ is non-decreasing and $\Prob$-a.s.~bounded, thus, it converges $\Prob$-a.s.~for all $t\in [0,\infty]$.\qed\\

To assess the type of convergence that we have for the sequence $Y^{v,k}$, we introduce a sequence of families of processes corresponding to a truncation of the time interval. For each $T>0$ and $k\geq 0$, we define the consistent family $((^T\!Y^{v,k}_t)_{t\geq 0}:v\in\mcU^f)$ of \cadlag supermartingales as
\begin{align*}
{}^T\!Y^{v,k}_t=\esssup_{u\in\mcU^k_{[t,T)}}\E\big[\varphi(v\circ u)-\sum_{j=1}^N c(v\circ[u]_{j})\big|\mcF_t\big]
\end{align*}
for all $v\in\mcU^f$ with $(t,b)\mapsto {}^T\!Y^{v\circ(t,b),k}_t\in\mcH'_\bbF$.
Then,

\begin{lem}\label{lem:TY}
The sequence $((^T\!Y^{v,k}_s)_{s\geq 0}:v\in\mcU^f)_{k\geq 0}$ satisfies:
\begin{enumerate}[i)]
  \item\label{lem-req:TY-cadlag} $Y^{v,0}\leq{^T}\!Y^{v,k}\leq Y^{v,k}$.
  \item For each $\epsilon>0$ there is a $T\geq 0$ such that $\Prob[\cup_{k=0}^\infty B^{T,\epsilon}_k]<\epsilon$, with
  \begin{align*}
    B^{T,\epsilon}_k:=\big\{\omega\in\Omega:\sup_{s\in[0,\infty]}\sup_{b\in U}|Y^{v\circ(s,b),k}_s-{^T}\!Y^{v\circ(s,b),k}_s|>\epsilon\big\}.
  \end{align*}
  \item There is a $\Prob$-a.s.~finite $\mcF$-measurable random variable $\xi$ and a constant $q>0$ such that
   \begin{align*}
     \sup_{t\in[0,\infty]}\sup_{b\in U}|^T\!Y^{v\circ (t,b),k}_t-{^T}\!Y^{v\circ(t,b),k'}_t|\leq \xi/(k')^{q}
   \end{align*}
   $\Prob$-a.s.~for each $0< k'\leq k$.
\end{enumerate}
\end{lem}

\noindent\emph{Proof.} The inequality in \emph{\ref{lem-req:TY-cadlag})} follows from noting that $\emptyset\in\mcU^k_{[t,T)}\subset\mcU^k_{t}$.\\

For the second statement we note that by Lemma~\ref{lem:mcH-prop}.(c), the process $(\sup_{b\in U}|Y^{v\circ(s,b),k}_s-{^T}\!Y^{v\circ(s,b),k}_s|:s\geq 0)$ is $\mcP_\bbF$-measurable and \cadlagp Now, each $u\in \mcU^k_{\tau}$ can be decomposed as $u=u_1\circ u_2$ with $u_1\in \mcU^k_{[\tau,T)}$ and $u_2\in \mcU^k_{T}$, which implies that
\begin{align}
Y^{v\circ(\tau,\beta),k}_\tau-{^T}\!Y^{v\circ(\tau,\beta),k}_\tau&=\esssup_{u\in\mcU^k_{\tau}}\E\big[\varphi(v\circ(\tau,\beta)\circ u)-\sum_{j=1}^N c(v\circ(\tau,\beta)\circ[u]_{j})\big|\mcF_\tau\big]\nonumber
\\
&\quad- \esssup_{u\in\mcU^k_{[\tau,T)}}\E\big[\varphi(v\circ(\tau,\beta)\circ u)-\sum_{j=1}^N c(v\circ(\tau,\beta)\circ[u]_{j})\big|\mcF_\tau\big]\nonumber
\\
&\leq \esssup_{u_1\in\mcU^k_{[\tau,T)},u_2\in\mcU^k_{T}}\E\big[\varphi(v\circ(\tau,\beta)\circ u_1\circ u_2)-\sum_{j=1}^{N_1+N_2} c(v\circ(\tau,\beta)\circ[u_1\circ u_2]_{j})\big|\mcF_\tau\big]\nonumber
\\
&\quad- \esssup_{u\in\mcU^k_{[\tau,T)}}\E\big[\varphi(v\circ(\tau,\beta)\circ u)-\sum_{j=1}^N c(v\circ(\tau,\beta)\circ[u]_{j})\big|\mcF_\tau\big]\nonumber
\\
&\leq \esssup_{u_1\in\mcU^k_{[\tau,T)},u_2\in\mcU^k_{T}}\E\big[|\varphi(v\circ(\tau,\beta)\circ u_1\circ u_2)-\varphi(v\circ(\tau,\beta)\circ u_1)|\big|\mcF_\tau\big],\label{ekv:TY-Y-bnd}
\end{align}
where $N_1$ and $N_2$ are the number of interventions in $u_1$ and $u_2$, respectively. We thus define the sets
\begin{align*}
\tilde B^{T,\epsilon}_k:=\Big\{\omega\in\Omega: \sup_{s\in[0,\infty]}\sup_{b\in U}\esssup_{u_1\in\mcU^k_{[0,T)},u_2\in\mcU^k_{T}}\E\big[|\varphi(v\circ(s,b)\circ u_1\circ u_2)-\varphi(v\circ(s,b)\circ u_1)|\big|\mcF_s\big]\geq\epsilon\Big\}
\end{align*}
and have by Equation \eqref{ekv:TY-Y-bnd} that $B^{T,\epsilon}_k\subset \tilde B^{T,\epsilon}_k$ for $k\geq 0$. Furthermore, as $\mcU^k_{[\tau,T)}\subset \mcU^{k+1}_{[\tau,T)}$ and $\mcU^k_{T}\subset \mcU^{k+1}_{T}$ we find that $\tilde B^{T,\epsilon}_k\subset \tilde B^{T,\epsilon}_{k+1}$ for all $k\geq 0$.

Now, let
\begin{align*}
\tau_k^\epsilon:=\inf\{s\geq 0:\sup_{b\in U}\esssup_{u_1\in\mcU^k_{[0,T)},u_2\in\mcU^k_{T}}\E\big[|\varphi(v\circ(s,b)\circ u_1\circ u_2)-\varphi(v\circ(s,b)\circ u_1)|\big|\mcF_s\big]\geq \epsilon\}
\end{align*}
(recalling our convention that $\inf\emptyset=\infty)$ and pick $\beta_k^\epsilon$ such that
\begin{align*}
  &\sup_{b\in U}\esssup_{u_1\in\mcU^k_{[0,T)},u_2\in\mcU^k_{T}}\E\big[|\varphi(v\circ(\tau_k^\epsilon,b)\circ u_1\circ u_2)-\varphi(v\circ(\tau_k^\epsilon,b)\circ u_1)|\big|\mcF_{\tau_k^\epsilon}\big]
  \\
  &\leq \esssup_{u_1\in\mcU^k_{[0,T)},u_2\in\mcU^k_{T}}\E\big[|\varphi(v\circ(\tau_k^\epsilon,\beta_k^\epsilon)\circ u_1\circ u_2)-\varphi(v\circ(\tau_k^\epsilon,\beta_k^\epsilon)\circ u_1)|\big|\mcF_{\tau_k^\epsilon}\big]+\epsilon/2.
\end{align*}
Then, by right continuity we have $\tilde B^{T,\epsilon}_k\subset \hat B^{T,\epsilon}_k$ where
\begin{align*}
\hat B^{T,\epsilon}_k:=\Big\{\omega\in\Omega: \esssup_{u_1\in\mcU^k_{[0,T)},u_2\in\mcU^k_{T}}\E\big[|\varphi(v\circ(\tau_k^\epsilon,\beta_k^\epsilon)\circ u_1\circ u_2)-\varphi(v\circ(\tau_k^\epsilon,\beta_k^\epsilon)\circ u_1)|\big|\mcF_{\tau_k^\epsilon}\big]\geq\epsilon/2\Big\}.
\end{align*}
We thus only need to show that there is a $T>0$ such that $\Prob[\hat B^{T,\epsilon}_k]<\epsilon$ for all $k\geq 0$. Arguing as in the proof of Proposition~\ref{prop:Ykbnd} we have
\begin{align*}
&\E\big[|\esssup_{u_1\in\mcU^k_{[0,T)},u_2\in\mcU^k_{T}}\E\big[|\varphi(v\circ(\tau_k^\epsilon,\beta_k^\epsilon)\circ u_1\circ u_2)-\varphi(v\circ(\tau_k^\epsilon,\beta_k^\epsilon)\circ u_1)|\big|\mcF_{\tau_k^\epsilon}\big]|^2\big]
\\
&\leq \E\big[\sup_{s\in[0,\infty]}|\esssup_{u_1\in\mcU_{[0,T)},u_2\in\mcU_{T}}\E\big[|\varphi(u_1\circ u_2)-\varphi(u_1)|\big|\mcF_{s}\big]|^2\big]
\\
&\leq C\sup_{u_1\in\mcU_{[0,T)},u_2\in\mcU_{T}}\E\big[|\varphi( u_1\circ u_2)-\varphi(u_1)|^2\big].
\end{align*}
We note that the right-hand side of the last inequality is independent of $k$ and by Definition~\ref{ass:onPSI}.(\ref{ass:@end}) it tends to 0 as $T\to\infty$. We thus conclude that there is a $T=T(\epsilon)$ such that $\Prob[\hat B^{T,\epsilon}_k]<\epsilon$ for all $k\geq 0$.\\

Concerning the third statement, we note that for $p\in(1,2)$, we have for each $\tau\in\mcT$ and all $\beta\in\mcA(\tau)$, that
\begin{align*}
{^T}\!Y^{v\circ(\tau,\beta),k}_\tau&\leq\sup_{s\in[0,\infty]}Y^{v\circ(\tau,\beta),k}_s
\\
&\leq \sup_{s\in[0,\infty]}\esssup_{u\in\mcU_s}\E[|\varphi(u)|\big|\mcF_s]
\\
&\leq 1+\sup_{s\in[0,\infty]}\esssup_{u\in\mcU_s}\E[|\varphi(u)|^p\big|\mcF_s]=: K(\omega)
\end{align*}
and similarly
\begin{align*}
{^T}\!Y^{v\circ(\tau,\beta),k}_\tau&\geq- K(\omega)
\end{align*}
for all $k\geq 0$ (where the inequalities hold $\Prob$-a.s.). Now, arguing as in the proof of Proposition~\ref{prop:Ykbnd} we have
\begin{align*}
&\E\big[\sup_{s\in[0,\infty]}\esssup_{u\in\mcU_s}\E[|\varphi(u)|^p\big|\mcF_s]^{2/p}\big]
\leq C\sup_{u\in\mcU}\E\big[|\varphi(u)|^2\big]<\infty
\end{align*}
and we conclude that there is a $\Prob$-null set $\mcN$ such that for each $\omega\in\Omega\setminus\mcN$ we have $K(\omega)<\infty$.

For $\epsilon>0$, let $u^{k,\epsilon}:=(\tau^{k,\epsilon}_1,\ldots,\tau^{k,\epsilon}_{N^{k,\epsilon}};\beta^{k,\epsilon}_1,\ldots,\beta^{k,\epsilon}_{N^{k,\epsilon}})\in \mcU^{k}_{[\tau,T)}$ be such that
\begin{align*}
{^T}\!Y^{v\circ(\tau,\beta),k}_\tau\leq\E\big[\varphi(v\circ(\tau,\beta)\circ u^{k,\epsilon})-\sum_{j=1}^{N^{k,\epsilon}} c(v\circ(\tau,\beta)\circ[u^{k,\epsilon}]_{j})\big|\mcF_\tau\big]+\epsilon.
\end{align*}
We note that on $[\tau\geq T]$ we have $N^{k,\epsilon}=0$ and get that for $\omega\in\Omega\setminus\mcN$ (in the remainder of the proof $\mcN$ denotes a generic $\Prob$-null set), we have
\begin{align*}
-K-\epsilon\leq \E\Big[\varphi(v\circ(\tau,\beta)\circ u^{k,\epsilon})-\sum_{j=1}^{N^{k,\epsilon}} c(v\circ(\tau,\beta)\circ[u^{k,\epsilon}]_{j})\Big|\mcF_\tau\Big]\leq K-\delta(T)\E\big[N^{k,\epsilon}\big|\mcF_\tau\big]
\end{align*}
and we conclude that $\E[\ett_{[N^{k,\epsilon}>k']}|\mcF_\tau]<(2K+\epsilon)/(\delta(T)k')$ for all $k'\geq 0$.

Now, for all $0\leq k'\leq k$ we have,
\begin{align*}
{^T}\!\breve Y^{v\circ(\tau,\beta),k,k'}_\tau&:=\E\Big[\varphi(v\circ(\tau,\beta)\circ [u^{k,\epsilon}]_{k'})-\sum_{j=1}^{N^{k,\epsilon}\wedge k'} c(v\circ(\tau,\beta)\circ[u^{k,\epsilon}]_{j})\big|\mcF_\tau\Big]
\\
&\leq {^T}\!Y^{v\circ(\tau,\beta),k'}_\tau\leq {^T}\!Y^{v\circ(s,b),k}_\tau,
\end{align*}
where we have introduced ${^T}\!\breve Y^{v\circ(\tau,\beta),k,k'}_\tau$ corresponding to the truncation $[u^{k,\epsilon}]_{k'}:=$\\$(\tau^{k,\epsilon}_1,\ldots,\tau^{k,\epsilon}_{N^{k,\epsilon}\wedge k'};\beta_1^{k,\epsilon},\ldots,\beta^{k,\epsilon}_{N^{k,\epsilon}\wedge k'})$ of $u^{k,\epsilon}$. As the truncation only affects the performance of the controller when $N^{k,\epsilon}>k'$ we have
\begin{align*}
{^T}\!Y^{v\circ(\tau,\beta),k}_\tau-{^T}\!\breve Y^{v\circ(\tau,\beta),k,k'}_\tau
&\leq \E\Big[\ett_{[N^{k,\epsilon}>k']}\big(\varphi(v\circ u^{k,\epsilon}) -\sum_{j=k'+1}^{N^{k,\epsilon}}c(v\circ(\tau,\beta)\circ[u^{k,\epsilon}]_{j})
\\
&\quad-\varphi(v\circ(\tau,\beta)\circ[u^{k,\epsilon}]_{k'})\big) \big|\mcF_\tau\Big]+\epsilon
\\
&\leq\E\Big[\ett_{[N^{k,\epsilon}>k']}\big(\varphi(v\circ(\tau,\beta)\circ u^{k,\epsilon})
-\varphi(v\circ(\tau,\beta)\circ[u^{k,\epsilon}]_{k'})\big) \big|\mcF_\tau\Big]+\epsilon.
\end{align*}
Applying H\"older's inequality we get that for $\omega\in\Omega\setminus\mcN$,
\begin{align*}
{^T}\!Y^{v\circ(\tau,\beta),k}_\tau-{^T}\!\breve Y^{v\circ(\tau,\beta),k,k'}_\tau&\leq 2\E[\ett_{[N^{k,\epsilon}>k']}|\mcF_\tau]^{1/q}
\esssup_{u\in\mcU}\E[|\varphi(v\circ(\tau,\beta)\circ u)|^p\big|\mcF_\tau]^{1/p}+\epsilon
\\
&\leq 2^{1+\frac{1}{q}}\frac{K(\omega)+\epsilon}{(\delta(T)k')^{1/q}}+\epsilon,
\end{align*}
with $\frac{1}{p}+\frac{1}{q}=1$. Since $\delta(T)>0$, there is thus a $\Prob$-a.s.~finite $\mcF$-measurable r.v.~$\xi=\xi(\omega)$ such that
\begin{align*}
|{^T}\!Y^{v\circ(\tau,\beta),k}_\tau-{^T}\!Y^{v\circ(\tau,\beta),k'}_\tau| \leq \frac{\xi}{(k')^{1/q}}+C\epsilon.
\end{align*}
Since, $\beta\in\mcA(\tau)$ was arbitrary we can choose $\beta$ such that
\begin{align*}
\sup_{b\in U}|{^T}\!Y^{v\circ(\tau,b),k}_\tau-{^T}\!Y^{v\circ(\tau,b),k'}_\tau|\leq |{^T}\!Y^{v\circ(\tau,\beta),k}_\tau-{^T}\!Y^{v\circ(\tau,\beta),k'}_\tau|+\epsilon \leq \frac{\xi}{(k')^{1/q}}+C\epsilon,
\end{align*}
$\Prob$-a.s.~and by right-continuity the last statement follows as $\epsilon>0$ was arbitrary.\qed\\

\begin{prop}\label{prop:conv}
For each $v\in\mcU^f$ we have
\begin{align*}
\sup_{t\in[0,\infty]}\sup_{b\in U} |Y^{v\circ(t,b),k}_t - \bar Y^{v\circ(t,b)}_t|\to 0
\end{align*}
as $k\to\infty$, outside of a $\Prob$-null set.
\end{prop}

\noindent\emph{Proof.} By Lemma~\ref{lem:TY}.(ii) there exist for each $\epsilon>0$ a $T\geq 0$ and a measurable set $B\subset\Omega$ with $\Prob[B]\geq 1-\epsilon$ such that
\begin{align*}
\sup_{t\in[0,\infty]}\sup_{b\in U} |Y^{v\circ(t,b),k}_t- Y^{v\circ(t,b),k'}_t|\leq \sup_{t\in[0,\infty]}\sup_{b\in U} |{^T}\!Y^{v\circ(t,b),k}_t-{^T}\!Y^{v\circ(t,b),k'}_t|+2\epsilon
\end{align*}
for all $0\leq k\leq k'$ and $\omega\in B$. Furthermore, by Lemma~\ref{lem:TY}.(iii) there is a $\Prob$-a.s.~finite r.v., $\xi$, such that
\begin{align*}
\sup_{t\in[0,\infty]}\sup_{b\in U} |Y^{v\circ(t,b),k}_t- Y^{v\circ(t,b),k'}_t|\leq \frac{\xi}{(k)^q}.
\end{align*}
Combining these and taking the limit as $k,k'\to\infty$ we find that
\begin{align*}
\lim_{k\to\infty}\sup_{t\in[0,\infty]}\sup_{b\in U} |\bar Y^{v\circ(t,b)}_t- Y^{v\circ(t,b),k}_t|\leq \lim_{k\to\infty}\lim_{k'\to\infty}\sup_{t\in[0,\infty]}\sup_{b\in U} |Y^{v\circ(t,b),k'}_t- Y^{v\circ(t,b),k}_t|\leq 2\epsilon
\end{align*}
on $B\setminus\mcN$ for some $\Prob$-null set $\mcN$. Now, as $\epsilon>0$ was arbitrary the statement follows.\qed\\

%
%

We are now ready to show that a verification family exists, establishing the existence of optimal controls for Problem 1.

\begin{prop}
A verification family exists.
\end{prop}

\noindent\emph{Proof.} Consistency of $((\bar Y^v_s)_{s\geq 0}:v\in\mcU^f)$ follows by consistency of each family in the sequence $((Y^{v,k}_s)_{s\geq 0}:v\in\mcU^f)_{k\geq 0}$ and uniform convergence. We treat each of the remaining properties in the definition of a verification family separately:\\

\noindent\emph{a)} By Proposition~\ref{prop:conv} and applying first (d) and then (b) of Lemma~\ref{lem:mcH-prop} we have that $\sup_{b\in U}\{- c_{s,b}^v+\bar Y^{v\circ(s,b)}_s\}$ is a \cadlagSTOP, quasi-left upper semi-continuous process of class [D]. In particular we note that
\begin{align*}
\ett_{[s= \infty]}\varphi(v)+\ett_{[s < \infty]}\sup_{b\in U}\{- c_{s,b}^v+\bar Y^{v\circ(s,b)}_s\}
\end{align*}
is \cadlagp Applying (iv) of Theorem \ref{thm:Snell} then gives
\begin{align*}
\bar Y^{v}_s=\esssup_{\tau \in \mcT_{s}} \E\Big[&\ett_{[\tau = \infty]}\varphi(v)+\ett_{[\tau < \infty]}\sup_{b\in U}\{- c_{s,b}^v+\bar Y^{v\circ(\tau,b)}_\tau\}\big| \mcF_s\Big].
\end{align*}

\bigskip

\noindent\emph{b)} Uniform boundedness was shown in Proposition~\ref{prop:Ykbnd}.

\bigskip

\noindent\emph{c)} As we noted in step a), this follows from Proposition~\ref{prop:conv} and Lemma~\ref{lem:mcH-prop}.(d).

\bigskip

\noindent\emph{d)} This is immediate from Proposition~\ref{prop:Ykdiscount}.\qed


\section{Application to impulse control of SFDEs\label{sec:SFDEs}}
In~\cite{SwitchElephant} a finite horizon impulse control problem with a discrete set $U$ was solved when the underlying process followed a stochastic delay differential equation (SDDE) under a particular loop condition on the impulses. This problem was motivated by hydropower operation where the flow-times between different power plants induce delays in the dynamics of the controlled system.

In this section we extend the results from~\cite{SwitchElephant} by considering a discounted infinite horizon setting, allowing an uncountable control set $U$ and also by taking the dynamics of the underlying process to follow a stochastic functional differential equation. Furthermore, our prior treatment of the problem with abstract reward, $\varphi$, and intervention cost, $c$, allows us to consider a less restrictive set of assumptions on the coefficients in the problem formulation. In particular we are able to remove the loop condition.

Our treatment of non-Markovian impulse control problems in infinite horizon should also be compared to~\cite{DjehicheInfHorImp} where an infinite horizon impulse control problem in a non-Markovian setting with a fixed discrete delay is considered. The work presented in this section goes in a different direction by having an underlying dynamics driven by a L\'evy process that is affected by the impulses in the control, resulting in a more complex relation between the control and the output of the performance functional. Furthermore, we investigate the important extension to random horizon which turns out to be a trivial modification of our initial problem.

Throughout this section, we will only consider controls for which $\tau_j\to\infty$, $\Prob$-a.s., and restrict our attention to the setting when the underlying uncertainty stems from a process $X^u$ with $u\in \mcU^f$ defined as $X^u:=\lim_{k\to\infty}X^{u,k}$ where
\begin{align}\label{ekv:SFDE1}
X^{u,0}_t&=x(t)\quad{\rm for}\:t\in (-\infty,0)
\\
X^{u,0}_t&=x(0)+\int_{0}^ta(s,(X^{u,0}_{r})_{r\leq s})ds +\int_{0}^t\sigma(s,(X^{u,0}_{r})_{r\leq s})dB_s\nonumber
\\
&\quad+\int_{0}^t\int_{\R^d\setminus\{0\}}\gamma(s,(X^{u,0}_{r})_{r< s},z)d\tilde P(ds,dz),\quad{\rm for}\:t\geq 0\label{ekv:SFDE2}
\end{align}
for some $x\in\bbD$, the set of all (deterministic) uniformly bounded, \cadlag functions $x:\R\to\R^d$, and
\begin{align}\label{ekv:SFDE3}
X^{u,k}_t&=X^{u,k-1}_t\quad{\rm for}\:t\in [0,\tau_k)
\\
X^{u,k}_t&=\Gamma(\tau_k,(X^{u,k-1}_{s})_{s\leq \tau_k},\beta_k)+\int_{\tau_k}^ta(s,(X^{u,k}_{r})_{r\leq s})ds +\int_{\tau_k}^t\sigma(s,(X^{u,k}_{r})_{r\leq s})dB_s\nonumber
\\
&\quad+\int_{\tau_k+}^t\int_{\R^d\setminus\{0\}}\gamma(s,(X^{u,k}_{r})_{r< s},z)d\tilde P(ds,dz),\quad{\rm for}\:t\geq \tau_k.\label{ekv:SFDE4}
\end{align}
The dynamics of $X$ are driven by a $d$-dimensional Brownian motion $B$ and a Poisson random measure $P$ with intensity measure $\varrho(ds; dz) = ds \times \mu(dz)$, where $\mu(dz)$ is the L\'evy measure on $\R^d$ of $P$ and $\tilde P(ds; dz) := (P - \varrho)(ds; dz)$ is called the compensated jump martingale random measure of $P$. We assume that $\bbG:=\{\mcG_t\}_{t\geq 0}$ is the natural filtration generated by $B$ and $P$, with $\mcG=\mcG_\infty:=\lim_{t\to\infty}\mcG_t$.

For notational simplicity, we assume that all uncertainty comes from the process $X$ and consider the discounted setting with a continuous discount factor $\rho:\R_+\to\R_+$. The reward functional is then of the type
\begin{align}\label{ekv:rewardSFDE}
J(u)=\E\bigg[\int_0^\infty e^{-\rho(t)}\phi(t,X^u_t)dt-\sum_{j=1}^Ne^{-\rho (\tau_j)}\ell(\tau_j,X^{u,j-1}_{\tau_j},\beta_j)\bigg].
\end{align}

\subsection{Assumptions}
We assume that the involved coefficients satisfy the following constraints:
\begin{ass}\label{ass:onSFDE}
For any $t,t'\geq 0$, $b,b'\in U$, $x,x'\in\R^d$ and $y,y'\in\bbD$ and for some $q\geq 2$ we have:
\begin{enumerate}[i)]
  \item\label{ass:onSFDE-Gamma} The function $\Gamma:\R_+\times\bbD\times U\to\R^d$ satisfies the Lipschitz condition
  \begin{align*}
    |\Gamma(t,(y_s)_{s\leq t},b)-\Gamma(t',(y'_s)_{s\leq t'},b')|&\leq C(\int_{-\infty}^{t\wedge t'}|y'_s-y_s|ds+|y'_{t'}-y_t|
    \\
    &\quad+(|t'-t|+|b'-b|)(\sup_{s\leq t}|y_s|+\sup_{s\leq t'}|y'_s|))
  \end{align*}
  and the growth condition
  \begin{align*}
    |\Gamma(t,(y_s)_{s\leq t},b)|\leq K_\Gamma\vee |y_t|.
  \end{align*}
  for some constant $K_\Gamma>0$.
  \item\label{ass:onSFDE-a-sigma} The coefficients $a:\R_+\times\bbD\to\R^{d}$ and $\sigma:\R_+\times\bbD\to\R^{d\times d}$ are continuous in $t$ and satisfy the growth condition
  \begin{align*}
    |a(t,(y_s)_{s\leq t})|+|\sigma(t,(y_s)_{s\leq t})|&\leq C(1+\sup_{s\leq t}|y_s|)
  \end{align*}
  and the Lipschitz continuity
  \begin{align*}
    |a(t,(y_s)_{s\leq t})-a(t,(y'_s)_{s\leq t})|+|\sigma(t,(y_s)_{s\leq t})-\sigma(t,(y'_s)_{s\leq t})|&\leq C\sup_{s\leq t}|y'_s-y_s|,
    \\
    \int_0^t |a(s,(y_r)_{r\leq s})-a(s,(y'_r)_{r\leq s})|ds&\leq C\int_{-\infty}^t|y'_s-y_s|ds
    \\
    \int_0^t |\sigma(s,(y_r)_{r\leq s})-\sigma(s,(y'_r)_{r\leq s})|^2ds&\leq C\int_{-\infty}^t|y'_s-y_s|^2ds.
  \end{align*}
  \item\label{ass:onSFDE-gamma} There is a $\bar \gamma:\R^d\to\R_+$, with $\int \bar \gamma^{2q}(z)\mu(dz)< \infty$ such that $\gamma:\R_+\times\bbD\times\R^d\to\R^{d}$ satisfies
  \begin{align*}
    |\gamma(t,(y_s)_{s< t},z)|&\leq \bar \gamma(z)(1+\sup_{s< t}|y_s|),
    \\
    |\gamma(t,(y_s)_{s< t},z)-\gamma(t,(y'_s)_{s< t}),z)|&\leq \bar \gamma(z)\sup_{s< t}|y_s-y_s'|,
    \\
    \int_0^t|\gamma(s,(y_r)_{r< s},z)-\gamma(s,(y'_r)_{r< s},z)|^{2(m+2)}ds &\leq \bar \gamma^{2(m+2)}(z)\int_{-\infty}^t|y'_s-y_s|^{2(m+2)}ds.
  \end{align*}
  \item\label{ass:onSFDE-phi} The running reward $\phi:\R_+\times \R^d\to\R$ is $\mcB(\R_+\times\R^d)$-measurable and satisfies the growth condition
  \begin{align*}
    |\phi(t,x)|\leq C(1+|x|^q)
  \end{align*}
  and for each $L>0$ there is a $C>0$ such that
  \begin{align*}
    |\phi(t,x)-\phi(t,x')|\leq C|x-x'|,
  \end{align*}
  whenever $|x|\vee|x'|\leq L$.
  \item\label{ass:onSFDE-c} There is a finite collection of closed connected subsets $(U_i)_{i= 1}^M$ of $U$ and corresponding maps $\ell_i:\R_+\times \R^d\times U_i\to \R$ that are jointly continuous in $(t,x,b)$, bounded from below, \ie
  \begin{align*}
    \ell_i(t,x,b)\geq\delta >0,
  \end{align*}
  of polynomial growth,
  \begin{align*}
    |\ell_i(t,x,b)|\leq C(1+|x|^q),
  \end{align*}
  and locally Lipschitz, \ie for each $L>0$ there is a $C>0$ such that
  \begin{align*}
    |\ell_i(t,x,b)-\ell_i(t,x',b)|\leq C|x-x'|,
  \end{align*}
  whenever $|x|\vee|x'|\leq L$, and we have $\ell(t,x,b)=\min_{i:b\in U_i}\ell_i(t,x,b)$.
\end{enumerate}
\end{ass}

\begin{rem}
To see that the above SFDE is a generalization of discrete delay SDDEs with Lipschitz coefficients note that if $\chi:\R_+\times (\R^{d})^k\to\R$ satisfies
\begin{align*}
|\chi(t,x_1,\ldots,x_{k+1})-\chi(t,x_1,\ldots,x_{k+1})|\leq C(|x_1-x_1'|+\cdots+|x_{k+1}-x_{k+1}'|)
\end{align*}
then for $l\geq 0$ we have
\begin{align*}
&\int_0^t|\chi(s,y_s,y_{s-\delta_1},\ldots,y_{s-\delta_k})-\chi(s,y'_s,y'_{s-\delta_1},\ldots,y'_{s-\delta_k})|^lds
\\
&\leq C\int_0^t(|y_s-y_s'|+|y_{s-\delta_1}-y_{s-\delta_1}'|\cdots+|y_{s-\delta_k}-y_{s-\delta_k}'|)^lds
\\
&\leq C\int_{-\infty}^t|y_s-y_s'|^lds.
\end{align*}
\end{rem}

\begin{rem}
In the above assumptions the involved coefficients are all deterministic. We remark that a trivial extension is to allow these to depend on $\omega$ as well in which case the coefficients in the Lipschitz conditions can be taken to be non-decreasing, $\Prob$-a.s.~finite, $\mcP_\bbG$-measurable \cadlag processes.
\end{rem}

The motivation for allowing intervention costs that are discontinuous in $b$ is the important application of production systems, where increasing the production beyond a certain threshold may necessitate a costly startup of additional production units.

\subsection{Existence of optimal controls}
In this section we show that the problem of maximizing the reward functional~\eqref{ekv:rewardSFDE} has a solution. Throughout we will, for notational simplicity, only consider the one-dimensional case ($d=1$), but we note that all results extend trivially to higher dimensions. We start with the following moment estimate:

\begin{prop}\label{prop:SFDEmoment}
Under Assumption~\ref{ass:onSFDE}, the SFDE \eqref{ekv:SFDE1}-\eqref{ekv:SFDE4} admits a unique solution for each $u\in\mcU$. Furthermore, the solution has moments of order $2q$ on compacts, in particular we have for $T>0$, that
\begin{align}\label{ekv:SFDEmoment}
\sup_{u\in\mcU^f}\E\big[\sup_{s\in[0,T]}|X^{u}_t|^{2q}\big]\leq C
\end{align}
and for each $v\in\mcU^f$, we have
\begin{align}\label{ekv:SFDEmoment2}
\sup_{u\in\mcU^f}\E\big[\sup_{(t,b)\in [0,T]\times U}|\E\big[\sup_{s\in[t,T]}|X^{v\circ (t,b)\circ u}_s|^{q}\big|\mcF_t\big]|^2\big]\leq C
\end{align}
where $C=C(T,q)$.
\end{prop}

\noindent\emph{Proof.} By repeated use of Theorem 3.2 in~\cite{AgramOksen} existence and uniqueness of solutions to \eqref{ekv:SFDE1}-\eqref{ekv:SFDE4} follows as $\tau_j\to\infty$, $\Prob$-a.s. 
By Assumption~\ref{ass:onSFDE}.(\ref{ass:onSFDE-Gamma}) we get, for $t\in [\tau_{j},T]$, using integration by parts, that
\begin{align*}
|X^{u,j}_t|^2&= |X^{u,j}_{\tau_{j}}|^2+2\int_{\tau_{j}+}^t X^{u,j}_{s-}dX^{u,j}_s+\int_{\tau_{j}+}^t d[X^{u,j},X^{u,j}]_s
\\
&\leq K^2_\Gamma\vee |X^{u,j-1}_{\tau_{j}}|^2+2\int_{\tau_{j}+}^t X^{u,j}_{s-} dX^{u,j}_s+\int_{\tau_{j}+}^t d[X^{u,j},X^{u,j}]_s.
\end{align*}
We note that if $|X^{u,j}_t|> K_\Gamma$ and $|X^{u,j}_s|\leq K_\Gamma$ for some $s\in[0,t)$ then there is a largest time $\zeta<t$ such that $|X^{u,j}_{\zeta-}|\leq K_\Gamma$. This means that during the interval $(\zeta,t]$ interventions will not increase the magnitude $|X^{u,j}|$. By induction, since $|x_{0}|$ is finite, we find that
\begin{align*}
|X^{u,j}_t|^2&\leq C+\sum_{i=0}^{j} \Big\{2\int_{\zeta\vee(\tilde\tau_{i}+)}^{t\wedge\tilde\tau_{i+1}} X^{u,i}_{s-}dX^{u,i}_s+\int_{\zeta\vee(\tilde\tau_{i}+)}^{t\wedge\tilde\tau_{i+1}} d[X^{u,i},X^{u,j}]_s\Big\}
\end{align*}
for all $t\in[0,T]$, where $\zeta=\sup\{s\geq 0 : |X^u_s|\leq K_\Gamma\}\vee 0$, $\tilde\tau_0+=0$, $\tilde\tau_i=\tau_i$ for $i=1,\ldots,j$ and $\tilde\tau_{j+1}=\infty$. Letting
\begin{align*}
R_t:=\sum_{i=0}^{j} \Big\{2\int_{\tilde\tau_{i}+}^{t\wedge\tilde\tau_{i+1}} X^{u,i}_{s-}d X^{u,i}_s+\int_{\tilde\tau_{i}+}^{t\wedge\tilde\tau_{i+1}} d[X^{u,i},X^{u,j}]_s\Big\}
\end{align*}
we thus find that
\begin{align*}
\E\Big[\sup_{s\in [t,T]}|X^{u,j}_s|^{2q}\Big|\mcF_t\Big]&\leq C +\E\Big[|X^{u,j}_t|^{2q}+2\sup_{s\in [t,T]}|R_s-R_t|^q\Big|\mcF_t\Big].
\end{align*}
Now, since $X^{u,i}$ and $X^{u,j}$ coincide on $[0,\tau_{i+1\wedge j+1})$ we have
\begin{align*}
\sum_{i=0}^{j}\int_{\tilde\tau_{i}+}^{t\wedge\tilde\tau_{i+1}} X^{u,i}_{s-} dX^{u,i}_s
&=\int_{0}^t X^{u,j}_{s}a(s,(X^{u,j}_r)_{r\leq s})ds+\int_{0}^{t}X^{u,j}_{s}\sigma(s,(X^{u,j}_r)_{r\leq s})dW_s\\
&\quad+\int_{0}^{t}\int_{\R^d\setminus\{0\}}X^{u,j}_{s-} \gamma(s,(X^{u,j}_r)_{r< s},z)\tilde P(ds,dz),
\end{align*}
and
\begin{align*}
\sum_{i=0}^{j} \int_{\tilde\tau_{i}+}^{t\wedge\tilde\tau_{i+1}} d[X^{u,i},X^{u,j}]_s&=\int_{0}^{t} \sigma^2(s,(X^{u,j}_r)_{r\leq s})ds\\
&\quad+\int_{0}^{t}\int_{\R^d\setminus\{0\}}\gamma^2(s,(X^{u,j}_r)_{r< s},z)P(ds,dz).
\end{align*}
From Assumption~\ref{ass:onSFDE}.(\ref{ass:onSFDE-a-sigma})-(\ref{ass:onSFDE-gamma}) and the Burkholder-Davis-Gundy inequality we get that
\begin{align*}
\E\Big[\sup_{t \leq s\leq T}|R_s-R_t|^q\big|\mcF_t\Big]
&\leq C\E\Big[|\int_{t}^T (1+\sup_{r\leq s}|X^{u,j}_{r}|^{4})ds|^{q/2}+|\int_{t}^T (1+\sup_{r\leq s}|X^{u,j}_{r}|^{2})ds|^{q}\big|\mcF_t\Big]
\\
&\leq C(1+ T^{q-1})\E\Big[ \int_{t}^T (1+\sup_{r\leq s}|X^{u,j}_{r}|^{2q})ds\big|\mcF_t\Big]
\end{align*}
and Gr\"onwall's lemma gives that
\begin{align}\label{ekv:moment_steg1}
\E\Big[\sup_{s\in[t,T]}|X^{u,j}_s|^{2q}\big|\mcF_t\Big]&\leq C(1+ \E\Big[\sup_{s\in[0,t]}|X^{u,j}_s|^{2q}\big|\mcF_t\Big]),
\end{align}
$\Prob$-a.s., where the constant $C=C(T,q)$ does not depend on $u$ or $j$ and \eqref{ekv:SFDEmoment} follows by letting $t=0$. Applying \eqref{ekv:moment_steg1} to the left-hand side of \eqref{ekv:SFDEmoment2} we get
\begin{align*}
\E\Big[\sup_{(t,b)\in [0,T]\times U}|\E\big[\sup_{s\in[t,T]}|X^{v\circ (t,b)\circ u}_s|^{q}\big|\mcF_t\big]|^2\Big]&\leq C(1+\E\Big[\sup_{t\in [0,T]}|\E\big[\sup_{s\in[0,t]}|X^{v}_s|^{q}\big|\mcF_t\big]|^2
\\
&\qquad+\sup_{(t,b)\in [0,T]\times U}|\E\big[|\Gamma(t,X^{v}_t,b)|^{q}\big|\mcF_t\big]|^2\Big])
\\
&\leq C(1+\E\Big[\sup_{t\in [0,T]}|\E\big[\sup_{s\in[0,t]}|X^{v}_s|^{q}\big|\mcF_t\big]|^2\Big])
\\
&\leq C(1+\E\big[\sup_{t\in [0,T]}|X^{v}_t|^{2q}\big])
\end{align*}
and the desired result follows from \eqref{ekv:SFDEmoment}.\qed\\

\begin{lem}\label{lem:SFDEflow}
For each $k\geq 0$, there is a $\Prob$-null set $\mcN$ such that for all $\omega\in\Omega\setminus\mcN$ and all $(\vect,\vecb)\in\mcD^k$ the limit $\lim_{(\vect',\vecb')\to (\vect,\vecb)}X^{\vect',\vecb'}$ exists in the topology of uniform convergence on compact subsets of $\R_+\setminus\{t_1,\ldots,t_k\}$. Furthermore, for all $(t,b)\in\R_+\times U$, we have
\begin{align*}
  \lim_{t'\searrow t}\lim_{b'\to b}\sup_{s\in[t',T]}|X^{v\circ(t',b')\circ u}_{s}-X^{v\circ(t,b)\circ u}_{s}|&=0,
\end{align*}
$\Prob$-a.s., for any $T\geq 0$, and $v\in\mcU^f$ and $u\in\mcU^k$ (with an exception set that is independent of $(t,b)$).
\end{lem}

\noindent\emph{Proof.} 
Our proof will rely on a pre-localization argument and we introduce the following non-decreasing sequence of stopping times
\begin{align*}
\kappa_K:=\inf\{s\geq 0: \int_{\R^d\setminus\{0\}}\bar\gamma(z)P(\{s\},dz)\geq K\},
\end{align*}
for $K\geq 0$ and set $\Lambda_K:=[0,\kappa_K)$. By, Assumption~\ref{ass:onSFDE}.(\ref{ass:onSFDE-gamma}) it then follows that $\kappa_K\to\infty$, $\Prob$-a.s.~as $K\to\infty$. Furthermore, we note that on $\Lambda_K$ the magnitude of the jumps of $X^u$ due to the Poisson jump integral of \eqref{ekv:SFDE1}-\eqref{ekv:SFDE4} are bounded by $C+K\sup_{s\leq t} |X^{u}_s|$ and repeating the argument in the proof of Proposition~\ref{prop:SFDEmoment} gives that
\begin{align}\label{ekv:all-moments}
\sup_{u\in\mcU^f}\E\Big[\sup_{s\in[0,T]\cap \Lambda_K}|X^{u}_t|^{l}\Big]\leq C,
\end{align}
for all $l\geq 0$.

For $0\leq t\leq t'\leq T$, we let ${}^{t,t'}\!X^{v}$ solve the SFDE \eqref{ekv:SFDE1}-\eqref{ekv:SFDE4} with integrand $(1-\ett_{(t,t']}(s))\gamma(s,\cdot,\cdot)$ in the jump part and let $j\geq 0$ be the largest integer such that $\tau_{j}\leq t'$. Then by Assumption~\ref{ass:onSFDE}.(\ref{ass:onSFDE-Gamma}) we have for $l=1,\ldots,j+1$ (recalling that $[u]_l=(\tau_1,\ldots,\tau_{N\wedge l};\beta_1,\ldots,\beta_{N\wedge l})$ is the truncation of $u$ limiting the number of interventions to $l$),
\begin{align*}
|{}^{t,t'}\!X^{v\circ[(t',b')\circ u]_l}_{t'}&-{}^{t,t'}\!X^{v\circ [(t,b)\circ u]_l}_{t'}| \leq C\Big(\int_{t}^{\tau_{l-1}}|{}^{t,t'}\!X^{v\circ[(t',b')\circ u]_{l-1}}_{s}-{}^{t,t'}\!X^{v\circ[(t,b)\circ u]_{l-1}}_{s}|ds
\\
&\quad+ |{}^{t,t'}\!X^{v\circ[(t',b')\circ u]_{l-1}}_{t'} -{}^{t,t'}\!X^{v\circ[(t,b)\circ u]_{l-1}}_{t'}| + |{}^{t,t'}\!X^{v\circ[(t,b)\circ u]_{l-1}}_{t'}-{}^{t,t'}\!X^{v\circ[(t,b)\circ u]_{l-1}}_{\tau_{l-1}}|
\\
&\quad+
((t'-\tau_{l-1})+\ett_{[l=1]}|b'-b|)(\sup_{s\leq t'}|{}^{t,t'}\!X^{v\circ[(t',b')\circ u]_{l-1}}_{s}|+\sup_{s\leq \tau_l}|{}^{t,t'}\!X^{v\circ[(t,b)\circ u]_{l-1}}_{s}|)\Big),
\end{align*}
with $\tau_0:=t$. We define ${}^1\!X^l:={}^{t,t'}\!X^{v\circ[(t,b)\circ u]_l}$, ${}^2\!X^l:={}^{t,t'}\!X^{v\circ[(t',b')\circ u]_l}$ and let $\delta  X^l:={}^2\!X^l-{}^1\!X^l$ and set $\delta X:=\delta  X^{k+1}$. Then, since the jump part is deactivated during $(t,t']$ and by \eqref{ekv:all-moments}, we have
\begin{align*}
\E\big[|\delta X_{t'}|^{2(m+2)}\big]&\leq C(|t'-t|^{m+2}+|b'-b|^{m+2}).
\end{align*}
For $l=j+2,\ldots,N$, we have by Assumption~\ref{ass:onSFDE}.(\ref{ass:onSFDE-Gamma}) that
\begin{align*}
|\delta X^{l+1}_{\tau_l}| &\leq C(\int_{0}^{\tau_l}|\delta X^{l}_{s}|ds + |\delta X^{l}_{\tau_l}|).
\end{align*}
Now, for $s\geq t'$,
\begin{align*}
\delta X_{s}&=\delta X_{t'}+\sum_{l=j}^{N}\int_{(t'\vee\tau_l)+}^{s\wedge\tau_{l+1}} d(\delta X^{l+1})_r+\sum_{l=j+1}^N\ett_{[s\geq\tau_l]}(\delta X^{l+1}_{\tau_l}-\delta X^{l}_{\tau_l}),
\end{align*}
with $\tau_{N+1}:=\infty$. Taking the absolute value on both sides we get
\begin{align*}
|\delta X_{s}|&\leq |\delta X_{t'}|+ C(\int_{t'}^s|a(r,({}^2\!X_\zeta)_{\zeta\leq r})-a(r,({}^1\!X_\zeta)_{\zeta\leq r})|dr
\\
&\quad +\sum_{l=j}^N |\int_{t'\vee\tau_l}^{\tau_{l+1}\wedge s}\sigma(r,({}^2\!X_\zeta)_{\zeta\leq r})-\sigma(r,({}^1\!X_\zeta)_{\zeta\leq r})dW_r|
\\
&\quad +\int_{t'+}^{s}|\gamma(r,({}^2\!X_\zeta)_{\zeta< r},z)-\gamma(r,({}^1\!X_\zeta)_{\zeta< s},z)|\tilde P(dr,dz)+\int_t^s|\delta X_{r}|dr).
\end{align*}
The Burkholder-Davis-Gundy inequality now gives
\begin{align*}
\E\Big[\sup_{r\in[t',s]}|\delta X_{r}|^{2(m+2)}\Big] &\leq C\E\Big[|\delta X_{t'}|^{2(m+2)}+ \Big(\int_{t'}^s|a(r,({}^2\!X_\zeta)_{\zeta\leq r})-a(r,({}^1\!X_\zeta)_{\zeta\leq r})|dr\Big)^{2(m+2)}
\\
& + \Big(\int_{t'}^{s}|\sigma(r,({}^2\!X_\zeta)_{\zeta\leq r})-\sigma(r,({}^1\!X_\zeta)_{\zeta\leq r})|^2dr\Big)^{m+2}
\\
& +\Big(\int_{t'+}^{s}|\gamma(r,({}^2\!X_\zeta)_{\zeta < r},z)-\gamma(r,({}^1\!X_\zeta)_{\zeta < r},z)|^2 \tilde P(dr,dz)\Big)^{m+2}+\int_t^s|\delta X_{r}|^{2(m+2)}dr\Big].
\end{align*}
Appealing to the boundedness of the jumps and the integral Lipschitz conditions on the coefficients then gives that
\begin{align*}
&\E\Big[\sup_{r\in[t',s]\cap\Lambda_K}|\delta X_r|^{2(m+2)}\Big]\leq C(\int_{t'}^{s}\E\Big[\sup_{\zeta\in[t',r]\cap\Lambda_K}|\delta X_\zeta|^{2(m+2)}\Big]dr+ |t'-t|^{m+2}+|b'-b|^{m+2}),
\end{align*}
for all $s\in[0,T]$. Now, Gr\"onwall's lemma gives
\begin{align*}
&\E\Big[\sup_{s\in[t',T]\cap\Lambda_K}|\delta X_s|^{2(m+2)}\Big] \leq C (|t'-t|^{m+2}+|b'-b|^{m+2}),
\end{align*}
where $C$ does not depend on $u\in\mcU^k$. Furthermore, for each $t\in[0,T]$ and each $\omega\in\Omega\setminus\mcN$ (for some $\Prob$-null set $\mcN$) there is a $t'>t$ such that $P(\omega,(t,t'],\R^d)=0$. Uniform convergence on $[t',T]\cap\Lambda_K$, thus, follows by applying a Kolmogorov continuity argument (see \eg Theorem 72 in Chapter IV of~\cite{Protter}) and uniform right-continuity follows as $\kappa_K\to\infty$, $\Prob$-a.s. 
The existence of limits follows similarly.\qed\\

%

\begin{defn}\label{def:trunc-reward}
For all $v,u\in\mcU^f$ we define the map $\Psi^{v,u}:\Omega\times\R_+\times U$ as
\begin{align*}
\Psi^{v,u}(t,b)&:=\int_0^Te^{-\rho(s)}\phi(s,X^{v\circ(t,b)\circ u}_s)ds-\sum_{j=1}^{N}e^{-\rho(\tau_j\vee t)}\ell(\tau_j\vee t,X^{v\circ(t,b)\circ [u]_{j-1}}_{\tau_j\vee t},\beta_j),
\end{align*}
 and for $T,L\geq 0$ we define the double truncation $\Psi^{v,u}_{T,L}:\Omega\times\R_+\times U$ of $\Psi^{v,u}$ as
\begin{align*}
\Psi^{v,u}_{T,L}(t,b)&:=\int_0^Te^{-\rho(s)}\phi(s,{}^L\!X^{v\circ(t,b)\circ u}_s)ds-\sum_{j=1}^{N(T-)}e^{-\rho(\tau_j\vee t)}\ell(\tau_j\vee t,{}^L\!X^{v\circ(t,b)\circ [u]_{j-1}}_{\tau_j\vee t},\beta_j),
\end{align*}
where $N(T-):=\max\{j:\tau_j<T\}$ and ${}^L\!X^u_s:=\frac{L}{L\vee |X^u_s|}X^u_s$.
\end{defn}


\begin{cor}\label{cor:SFDEvarphi-cont}
For each $T,L\geq 0$, $k\geq 0$, $v\in\mcU^f$ and $u\in\mcU^k$ the map $(t,b)\mapsto\Psi^{v,u}_{T,L}(t,b)$ has limits everywhere and is $\Prob$-a.s.~continuous on $[\eta_j,\eta_{j+1})\times U$ where $(\eta_j)_{j\geq 1}$ are the jump times of $P$.
\end{cor}

\noindent\emph{Proof.} Let
\begin{align*}
\varphi_{T,L}(u)&:=\int_0^Te^{-\rho(s)}\phi(s,{}^L\!X^u_s)ds
\end{align*}
and note that for $0\leq t\leq t'$ we have
\begin{align*}
|\varphi_{T,L}(v\circ(t',b')\circ u)-\varphi_{T,L}((v\circ(t,b)\circ u)|&\leq \int_t^{t'}e^{-\rho(s)}(|\phi(s,{}^L\!X^{v\circ(t,b)\circ u}_s)|+|\phi(s,{}^L\!X^{v\circ(t',b')\circ u}_s)|)ds
\\
&\quad + \int_{t'}^{T}e^{-\rho(s)}(|\phi(s,{}^L\!X^{v\circ(t,b)\circ u}_s)-\phi(s,{}^L\!X^{v\circ(t',b')\circ u}_s)|)ds
\\
&\leq C(|t'-t|+\int_{t'}^{T}e^{-\rho(s)}(|{}^L\!X^{v\circ(t,b)\circ u}_s-{}^L\!X^{v\circ(t',b')\circ u}_s)|ds)
\\
&\leq C(|t'-t|+\int_{t'}^{T}e^{-\rho(s)}(|X^{v\circ(t,b)\circ u}_s-X^{v\circ(t',b')\circ u}_s)|ds).
\end{align*}
Now, by Lemma~\ref{lem:SFDEflow} it follows immediately that
\begin{align*}
\lim_{t'\searrow t}\lim_{b'\to b}|\varphi_{T}(v\circ(t',b')\circ u)-\varphi_{T}(v\circ(t,b)\circ u)|=0,
\end{align*}
and from its proof we have that
\begin{align*}
\lim_{t'\nearrow t}\lim_{b'\to b}|\varphi_{T}(v\circ(t',b')\circ u)-\varphi_{T}(v\circ(t,b)\circ u)|=0,
\end{align*}
whenever $t\notin\{\eta_1,\eta_2,\ldots\}$. Concerning the intervention costs we have
\begin{align}\nonumber
&\sum_{j=1}^{N(T-)} |e^{-\rho(\tau_j\vee t')}\ell(\tau_j\vee t',{}^L\!X^{v\circ(t',b')\circ [u]_{j-1}}_{\tau_j\vee t'},\beta_j)-e^{-\rho(\tau_j\vee t)}\ell(\tau_j\vee t,{}^L\!X^{v\circ(t,b)\circ [u]_{j-1}}_{\tau_j\vee t},\beta_j)|
\\
&\leq \sum_{j=1}^{N(T-)} \{\ett_{[0,t')}(\tau_j)|e^{-\rho( t')}\ell(t',{}^L\!X^{v\circ (t,b)\circ[u]_{j-1}}_{t'},\beta_j)-e^{-\rho(\tau_j\vee t)}\ell(\tau_j\vee t,{}^L\!X^{v\circ (t,b)\circ[u]_{j-1}}_{\tau_j\vee t},\beta_j)|\nonumber
\\
&\qquad+|\ell(\tau_j\vee t',{}^L\!X^{v\circ (t',b')\circ[u]_{j-1}}_{\tau_j\vee t'},\beta_j)-\ell(\tau_j\vee t',{}^L\!X^{v\circ (t,b)\circ[u]_{j-1}}_{\tau_j\vee t'},\beta_j)|\}\label{ekv:c-updelad}
\end{align}
where the first term tends to zero as $t'\searrow t$ by joint continuity of $\ell$, continuity of $\rho$ and right continuity of $X$. By continuity of $\ell$, the assertion follows by repeating the argument in the proof of Lemma~\ref{lem:SFDEflow}.\qed\\

\bigskip

\begin{lem}\label{lem:J-star-T}
For each $T>0$ and $k\geq 0$ there is, for every $v\in\mcU^f$, a $J^*_T\in\mcH_\bbF$ such that for all $\tau\in\mcT$ and $b\in U$ we have
\begin{align*}
J^*_T(\tau,b)=\esssup_{u\in\mcU^k}\E\Big[\Psi^{v,u}_T(\tau,b)\big|\mcF_\tau\Big],
\end{align*}
$\Prob$-a.s.~(with an exception set that is independent of $b$).
\end{lem}

\noindent\emph{Proof.} For any $K,L\geq 0$ it follows by Corollary~\ref{cor:SFDEvarphi-cont} and Theorem~\ref{thm:OP} that there is for each $b\in U$ an $\bbF$-optional \cadlag process $(Z^{b,u}_t:t\geq 0)$ such that
\begin{align*}
Z^{b,u}_\tau=\E\Big[\Psi^{v,u}_{T\wedge\kappa_K,L}(\tau,b)\big|\mcF_\tau\Big],
\end{align*}
$\Prob$-a.s. Now, pick a sequence $(\epsilon_l)_{l\geq 0}$ of positive real numbers such that $\epsilon_l\searrow 0$ and for $j,l\geq 0$ define $s^l_j:=j2^{-l}$. Then, there is a control $u^l_j\in \mcU^k$ such that
\begin{align*}
\E\Big[\Psi^{v,u^l_j}_{T\wedge\kappa_K,L}(\tau,b)\big|\mcF_\tau\Big]\geq\esssup_{u\in\mcU^k}\E\Big[\Psi^{v,u}_{T\wedge\kappa_K,L}(\tau,b)\big|\mcF_\tau\Big]-\epsilon_l.
\end{align*}
Define the sequence of \cadlag processes $(\tilde Z^{b,l}_t:t\geq 0)_{l\geq 0}$ as
\begin{align*}
\tilde Z^{b,l}_t:=\sum_{j=0}^{\infty}\ett_{[s^l_j,s^l_{j+1})}(t)Z^{b,u^l_{j+1}}_t
\end{align*}
and set $\hat Z^{b,l}_t:=\max_{i\in \{0,\ldots,l\}}\tilde Z^{b,i}_t$. Then, $\hat Z^{b,l}$ is an increasing $\Prob$-a.s.~finite sequence of \cadlag processes and it, thus, converges pointwisely, $\Prob$-a.s.~to a limit $Z^{b,*}$ that, moreover, is $\mcP_\bbF$-measurable. We note that for any $l\geq 0$ and $\tau\in\mcT^f$ we have with $\tau_l:=\inf\{s\geq\tau:s\in\Pi_l\}$ and $u^l:=\sum_{j=1}^\infty\ett_{[\tau_l=s^l_j]}u^l_j$, that
\begin{align*}
\esssup_{u\in\mcU^k}\E\Big[\Psi^{v,u}_{T\wedge\kappa_K,L}(\tau,b)\big|\mcF_\tau\Big] -Z^{b,*}_\tau &\leq\esssup_{u\in\mcU^k}\E\Big[\Psi^{v,u}_{T\wedge\kappa_K,L}(\tau,b)\big|\mcF_\tau\Big] -\E\Big[\Psi^{v,u^l}_{T\wedge\kappa_K,L}(\tau,b)\big|\mcF_\tau\Big]
\\
&=\esssup_{u\in\mcU^k}\E\Big[\Psi^{v,u}_{T\wedge\kappa_K,L}(\tau,b)\big|\mcF_\tau\Big] - \esssup_{u\in\mcU^k}\E\Big[\Psi^{v,u}_{T\wedge\kappa_K,L}(\tau_l,b)\big|\mcF_\tau\Big]
\\
&\quad + \esssup_{u\in\mcU^k}\E\Big[\Psi^{v,u}_{T\wedge\kappa_K,L}(\tau_l,b)\big|\mcF_\tau\Big] - \E\Big[\Psi^{v,u^l}_{T\wedge\kappa_K,L}(\tau,b)\big|\mcF_\tau\Big]
\end{align*}
and as
\begin{align*}
&|\esssup_{u\in\mcU^k}\E\Big[\Psi^{v,u}_{T\wedge\kappa_K,L}(\tau_l,b)\big|\mcF_\tau\Big] - \E\Big[\Psi^{v,u^l}_{T\wedge\kappa_K,L}(\tau,b)\big|\mcF_\tau\Big]|
\\
&\leq \esssup_{u\in\mcU^k}\E\Big[|\Psi^{v,u}_{T\wedge\kappa_K,L}(\tau,b)-\Psi^{v,u}_{T\wedge\kappa_K,L}(\tau_l,b)|\big|\mcF_\tau\Big] + \epsilon_l
\end{align*}
we get that
\begin{align*}
|Z^{b,*}_\tau-\esssup_{u\in\mcU^k}\E\Big[\Psi^{v,u}_{T\wedge\kappa_K,L}(\tau,b)\big|\mcF_\tau\Big]|&\leq 2\esssup_{u\in\mcU^k}\E\Big[|\Psi^{v,u}_{T\wedge\kappa_K,L}(\tau,b)-\Psi^{v,u}_{T\wedge\kappa_K,L}(\tau_l,b)|\big|\mcF_\tau\Big] + \epsilon_l.
\end{align*}
Doob's maximal inequality now gives that
\begin{align*}
\E\big[|Z^{b,*}_\tau-\esssup_{u\in\mcU^k}\E\Big[\Psi^{v,u}_{T\wedge\kappa_K,L}(\tau,b)\big|\mcF_\tau\Big]|^2\big]&
\leq C\sup_{u\in\mcU^k}\E\big[|\Psi^{v,u}_{T\wedge\kappa_K,L}(\tau,b)-\Psi^{v,u}_{T\wedge\kappa_K,L}(\tau_l,b)|^2\big]+\epsilon_l.
\end{align*}
Letting $l\to\infty$ and using that the map $t\mapsto \Psi^{v,u}_{T\wedge\kappa_K,L}(t,b)$ is right-continuous uniformly in $u$ it follows that
\begin{align*}
  Z^{b,*}_\tau=\esssup_{u\in\mcU^k}\E\Big[\Psi^{v,u}_{T\wedge\kappa_K,L}(\tau,b)\big|\mcF_\tau\Big],
\end{align*}
$\Prob$-a.s., for each stopping time $\tau\in\mcT$.

We now show that $Z^{b,*}$ is a \cadlag process. First, since $Z^{b,*}$ is the limit of an increasing sequence of \cadlag processes we have that $\liminf_{t'\searrow t}Z^{b,*}_{t'}\geq Z^{b,*}_{t}$. For any $\tau\in\mcT^f$ and $\epsilon>0$ let
\begin{align*}
\tau_l:=\inf\{s\geq\tau:\hat Z^{b,l}_s\geq Z^{b,*}_\tau+\epsilon\}.
\end{align*}
Then as $(Z^{b,l})_{l\geq 0}$ is non-decreasing, the sequence $(\tau_l)_{l\geq 0}$ is non-increasing. Let $B:=\{\omega\in\Omega:\lim_{l\to\infty}\tau_l=\tau\}$ and note that $B\in\mcF_\tau$ by right-continuity of the filtration and $\limsup_{t'\searrow t}Z^{b,*}_{t'}< Z^{b,*}_{t}+\epsilon$ on $B^c:=\Omega\setminus B$. Moreover, with $\hat\tau_l:=\ett_{B} \tau_l+\ett_{B^c}\tau$, Fatou's lemma gives
\begin{align*}
\liminf_{l\to\infty}\E\Big[Z^{b,*}_{\hat\tau_l}-Z^{b,*}_{\tau}\Big] = \liminf_{l\to\infty}\E\Big[\ett_{B}(Z^{b,*}_{\hat\tau_l}-Z^{b,*}_{\tau})\Big]\geq \Prob[B]\epsilon.
\end{align*}
On the other hand, we have
\begin{align*}
\E\Big[Z^{b,*}_{\hat\tau_l}-Z^{b,*}_{\tau}\Big]\leq\sup_{u\in\mcU^k}\E\big[\Psi^{v,u}_{T\wedge\kappa_K,L}(\tau,b) - \Psi^{v,u}_{T\wedge\kappa_K,L}(\hat\tau_l,b)|\big]\to 0
\end{align*}
as $l\to\infty$ and we conclude that $\Prob[B]=0$ and, since $\epsilon>0$ was arbitrary, it follows that $\limsup_{t'\searrow t}Z^{b,*}_{t'}= Z^{b,*}_{t}$.

To prove that $Z^{b,*}$ has left limits we define, for $\epsilon>0$, the sequence $(\vartheta^\epsilon_j)_{j\geq 0}$ as $\vartheta^\epsilon_0=0$ and then recursively let
\begin{align*}
\vartheta^\epsilon_j:=\inf\{s\geq \vartheta^\epsilon_{j-1}:\esssup_{u\in\mcU^k}|Z^{b,u}_s-Z^{b,u}_{\vartheta^\epsilon_{j-1}}|\geq \epsilon\}.
\end{align*}
We note by the above discussion that $\vartheta^\epsilon_j\in\mcT$ and furthermore, by right-continuity that $\vartheta^\epsilon_{j}> \vartheta^\epsilon_{j-1}$ and $\vartheta^\epsilon_j\to\infty$, $\Prob$-a.s. If not, we would have $\vartheta^\epsilon_j\nearrow \vartheta^\epsilon\in\mcT^f$ on some set $A$ of positive measure. However, as increments in the jump integral part is $\Prob$-a.s.~zero at predictable times we note by Corollary~\ref{cor:SFDEvarphi-cont} that $\Psi^{v,u}(t,b)$ is continuous in $t$ at $\vartheta^\epsilon$ on $A\setminus\mcN$ for some $\Prob$-null set $\mcN$, uniformly in $u$. Now, as the filtration is quasi-left continuous this implies that
\begin{align*}
\limsup_{j\to\infty}\esssup_{u\in\mcU^k}\sup_{b\in U}|Z^{b,u}_{\vartheta^\epsilon_{j}}-Z^{b,u}_{\vartheta^\epsilon_{j-1}}|=0,
\end{align*}
on $A\setminus\mcN$, a contradiction. Letting,
\begin{align*}
\breve{Z}^{b,l}_t:=\sum_{j=0}^{\infty}\ett_{[\vartheta^{1/l}_{j},\vartheta_{j-1})}(t)Z^{b,*}_{\vartheta^{1/l}_{j}},
\end{align*}
we find that $(\breve{Z}^{b,l})_{l\geq 0}$ is a sequence of \cadlag processes with $\sup_{t\in[0,T]}|Z^{b,*}_t-\breve{Z}^{b,l}_t|\leq 1/l$ and we conclude that $Z^{b,*}$ is \cadlagp

By repeating the argument in the proof of Lemma~\ref{lem:SFDEflow} we find that
\begin{align*}
\sup_{u\in\mcU^k}\E\Big[|\Psi^{v,u}_{T\wedge\kappa_K,L}(t,b')-\Psi^{v,u}_{T\wedge\kappa_K,L}(t,b)|^{m+1}\Big]\leq C|b'-b|^{m+1},
\end{align*}
$\Prob$-a.s.~for any $t\in\R_+$ and $b,b'\in U$ and it follows that
\begin{align*}
\E\Big[\sup_{t\in[0,\infty]}|Z^{b',*}_t-Z^{b,*}_t|^{m+1}\Big]\leq C|b'-b|^{m+1}.
\end{align*}
Hence, by Kolmogorov's continuity theorem and Lemma~\ref{cor:SFDEvarphi-cont} it follows that there is a unique map $h^{K,L}\in\mcH_\bbF$ such that
\begin{align*}
h^{K,L}(\tau,b)=\esssup_{u\in\mcU^k}\E\Big[\Psi^{v,u}_{T\wedge\kappa_K,L}(\tau,b)\big|\mcF_\tau\Big],
\end{align*}
$\Prob$-a.s. By dominated convergence we find that $h^{K,L}$ converges pointwisely to some $h$ as $K,L\to\infty$. We define the set
\begin{align*}
\Xi_L:=\{\omega\in\Omega:\sup_{(t,b)\in[0,T]\times U}\E\big[\sup_{s\in[t,T]}|X^{v\circ (t,b)\circ u}_s|\big|\mcF_t\big]>L\}
\end{align*}
and note that for $p\in (1,2)$, we have
\begin{align*}
&\E\Big[\sup_{(t,b)\in [0,T]\times U}|\esssup_{u\in\mcU^k}\E\Big[\Psi^{v,u}_{T}(t,b)\big|\mcF_t\Big]-h^{K,L}(t,b)|^p\Big]
\\
&\leq \sup_{u\in\mcU^k}\E\Big[\sup_{(t,b)\in [0,T]\times U}|\E\big[\Psi^{v,u}_{T}(t,b)-\Psi^{v,u}_{T\wedge\kappa_K,L}(t,b)\big|\mcF_t\big]|^p\Big]
\\
&\leq C\sup_{u\in\mcU^k}\E\big[(\ett_{[\kappa_K<T]}+\ett_{\Xi})\sup_{(t,b)\in [0,T]\times U}|\E\big[\sup_{s\in[t,T]}|X^{v\circ (t,b)\circ u}_s|^{q}\big|\mcF_t\big]|^p\big]
\\
&\leq C (\E\big[\ett_{[\kappa_K<T]}]^{1/p'}+\sup_{u\in\mcU^k}\E\big[\ett_{\Xi_L}]^{1/p'})
\end{align*}
where $\frac{1}{p'}+\frac{p}{2}=1$ and the last step follows by H\"older's inequality and Proposition~\ref{prop:SFDEmoment}. Now, the right-hand side of the last inequality goes to zero as $K,L\to\infty$ by the definition of $\kappa_K$ and Proposition~\ref{prop:SFDEmoment} and by uniform convergence we conclude that there is a $J^*_T\in\mcH_\bbF$ such that
\begin{align*}
J^*_T(\tau,b)=\esssup_{u\in\mcU^k}\E\Big[\Psi^{v,u}_{T}(\tau,b)\big|\mcF_\tau\Big],
\end{align*}
$\Prob$-a.s.

It remains to show that we can choose the exception set to be independent of $b$. Let $\bar U_0\subset \bar U_1\subset\cdots$ be a sequence of finite subsets of $U$ with $\min_{b\in \bar U_l}\max_{b'\in U}|b'-b|\leq 2^{-l}$. For $\beta\in\mcA(\tau)$ define $(\beta_l)_{l\geq 0}$ as a measurable selection of $\beta_l\in\argmin_{b\in\bar U_l}|\beta-b|$. Then since $\beta_l$ takes values in a finite set we have
\begin{align*}
J^*_T(\tau,\beta_l)=\esssup_{u\in\mcU^k}\E\Big[\Psi^{v,u}_{T}(\tau,\beta_l)\big|\mcF_\tau\Big],
\end{align*}
$\Prob$-a.s. By continuity it follows that
\begin{align*}
\lim_{l\to\infty}J^*_T(\tau,\beta_l)=J^*_T(\tau,\beta),
\end{align*}
$\Prob$-a.s. Furthermore, by uniform integrability and $\Prob$-a.s.~continuity of $\Psi^{v,u}_{T}$ uniformly in $u$ we have that
\begin{align*}
\lim_{l\to\infty}\esssup_{u\in\mcU^k}\E\big[|\Psi^{v,u}_{T}(\tau,\beta)-\Psi^{v,u}_{T}(\tau,\beta_l)|\big|\mcF_\tau\big]=0
\end{align*}
and we conclude that
\begin{align*}
J^*_T(\tau,\beta)=\esssup_{u\in\mcU^k}\E\Big[\Psi^{v,u}_{T}(\tau,\beta)\big|\mcF_\tau\Big],
\end{align*}
$\Prob$-a.s. From this the statement follows as $\beta\in\mcA(\tau)$ was arbitrary.\qed\\

This far we have not made any assumption on the discount factor $\rho$, other than it being continuous. Clearly, some assumptions on the growth of $\rho$ have to be made in order for the maximization problem to have a finite value. We summarize this in the following hypothesis:

\begin{hyp*}\textbf{[Disc.-A]}
There is an $\epsilon>0$ such that $\rho(t)\geq \epsilon t$,
\begin{align*}
\sup_{u\in\mcU^f}\E\Big[\int_T^\infty e^{-\rho(t)}|\phi(t,X^u_t)|^2dt\Big]\leq Ce^{-\epsilon T}
\end{align*}
and
\begin{align*}
\sup_{u\in\mcU^f}\E\Big[\sup_{t\in [0,\infty)}e^{-2\rho(t)}|\ell(t,X^u_t,b)|^2dt\Big]<\infty
\end{align*}
for all $T\geq 0$ and $b\in U$. Furthermore, for each $k\geq 0$ there is an $\epsilon>0$ such that for all $T\geq T'$ for some $T'>0$ we have
\begin{align*}
\sup_{u\in\mcU^k}\E\Big[\sup_{(t,b)\in [0,\infty)\times U}\E\Big[|\int_T^\infty e^{-\rho(s)}\phi(s,X^{v\circ(t,b)\circ u}_s)ds|\Big|\mcF_t\Big]\Big]\leq e^{-\epsilon T},
\end{align*}
and
\begin{align*}
\sup_{u\in\mcU^k}\E\Big[\sup_{(t,b)\in [0,\infty)\times U}\E\Big[\ett_{[N\geq 1]}\ett_{[\tau_N\vee t\geq T]}e^{-\rho(\tau_N\vee t)}|\ell(\tau_N,X^{v\circ(t,b)\circ u}_{\tau_N\vee t},\beta_N)|\Big|\mcF_t\Big]\Big]\leq e^{-\epsilon T}.
\end{align*}
for all $v\in\mcU^f$.
\end{hyp*}

\bigskip

\begin{rem}
An important situation where Hypothesis~\textbf{Disc.-A} holds with $\rho(t)=\rho_0 t$ for any $\rho_0>0$ is when the functions $\phi$ and $\ell$ are eventually bounded, \ie when there is a $T'>0$ such that $|\phi(t,x)|\leq C$ and $|\ell(u\circ(t,x))|\leq C$ for all $(t,x)\in [T',\infty)\times \R^d$. Another important case is when $\rho(T)-\ln(C(T,q))$ grows linearly in $T$, where $C(T,q)$ is the growth in Proposition~\ref{prop:SFDEmoment}.\\
\end{rem}

We are now ready to state the main result of this section, showing that under Assumption~\ref{ass:onSFDE} and Hypothesis~\textbf{Disc.-A} an optimal control for the problem of maximizing $J$ exist.

\begin{prop}\label{prop:SFDE-OC}
Under Hypothesis~\textbf{Disc.-A} there is a $u^*\in\mcU^f$ such that $J(u^*)\geq J(u)$ for all $u\in\mcU^f$. Furthermore, $u^*$ is given by the recursion \eqref{ekv:taujDEF}-\eqref{ekv:betajDEF}, with
\begin{align*}
\varphi(u):=\int_0^\infty e^{-\rho(t)}\phi(t,X^u_t)dt
\end{align*}
and
\begin{align*}
c(u\circ(t,b)):=e^{-\rho(t\vee\tau_N)}\ell(t\vee\tau_N,X^u_{t\vee\tau_N},b).
\end{align*}
\end{prop}

\noindent\emph{Proof.} To show that the assertion is true we need to show that the pair $(\varphi,c)$ is an admissible reward pair. It is clear that the uniform $L^2$-bounds on $\varphi$ and $c$ in Definition~\ref{ass:onPSI}.(\ref{ass:psiBND}) hold by Hypothesis~\textbf{Disc.-A}. In particular, we note that by Jensen's inequality we get
\begin{align*}
\E[|\varphi(u)|^2]&=\E\Big[|\int_0^\infty e^{-\rho(t)}\phi(t,X^u_t)dt|^2\Big]
\\
&\leq C\E\Big[\int_{0}^{\infty} e^{-\rho(t)}|\phi(t,X^u_t)dt|^2dt\Big]\leq C.
\end{align*}

The decreasing importance property stated in Definition~\ref{ass:onPSI}.(\ref{ass:@end}) follows similarly by noting that for $v\in\mcU^f_T$ with $T\geq T'$ we have, by Hypothesis~\textbf{Disc.-A}, that
\begin{align*}
\E[|\varphi(u\circ v)-\varphi(u)|^2]&=\E\Big[|\int_T^\infty e^{-\rho(t)}(\phi(t,X^{u\circ v}_t)-\phi(t,X^u_t))dt|^2\Big]
\\
&\leq C\E\Big[\int_T^\infty e^{-\rho(t)}(|\phi(t,X^{u\circ v}_t)|^2+|\phi(t,X^u_t)|^2) dt\Big]
\\
&\leq Ce^{-\epsilon T},
\end{align*}
which tends to 0 as $T\to\infty$.

Concerning the continuity properties listed in Definition~\ref{ass:onPSI}.(\ref{ass:psiREG}) we note that for each $k\geq 0$ and $v\in\mcU^f$ we have
that
\begin{align*}
|\Psi^{v,u}(t',b')-\Psi^{v,u}(t,b)|&\leq |\Psi^{v,u}_T(t',b')-\Psi_T^{v,u}(t,b)|+|\Psi^{v,u}(t',b')-\Psi_T^{v,u}(t',b')|
\\
&\quad+|\Psi^{v,u}(t,b)-\Psi_T^{v,u}(t,b)|.
\end{align*}
Now,
\begin{align*}
&\sup_{u\in\mcU^k}\E\Big[\sup_{(t,b)\in\R_+\times U}\E\Big[|\varphi(v\circ(t,b)\circ u)-\varphi_T(v\circ(t,b)\circ u)|\Big|\mcF_t\Big]\Big]
\\
&=\sup_{u\in\mcU^k}\E\Big[\sup_{(t,b)\in \R_+\times U}\E\Big[|\int_T^\infty e^{-\rho(s)}\phi(s,X^{v\circ(t,b)\circ u}_s)ds|\Big|\mcF_t\Big]\Big]\leq e^{-\epsilon T}.
\end{align*}
and similarly
\begin{align*}
&\sup_{u\in\mcU^k}\E\Big[\sup_{(t,b)\in\R_+\times U}\E\Big[|c(v\circ(t,b)\circ u)-\ett_{[\tau_N\vee t < T]}c(v\circ(t,b)\circ u)|\Big|\mcF_t\Big]\Big]
\\
&=\sup_{u\in\mcU^k}\E\Big[\sup_{(t,b)\in \R_+\times U}\E\Big[\ett_{[N\geq 1]}\ett_{[\tau_N\vee t\geq T]}e^{-\rho(\tau_N\vee t)}|\ell(\tau_N\vee t,X^{v\circ(t,b)\circ u}_{\tau_N\vee t},\beta_N)|\Big|\mcF_t\Big]\Big]\leq e^{-\epsilon T}.
\end{align*}
This implies that
\begin{align*}
&\Prob\Big[\sup_{(t,b)\in\R_+\times U}\E\big[|\Psi^{v,u}(t,b)-\Psi^{v,u}_T(t,b)|\big|\mcF_t\big]\geq e^{-\epsilon T/2}\Big]\leq Ce^{-\epsilon T/2}
\end{align*}
and the Borel-Cantelli lemma gives that $\sup_{(t,b)\in\R_+\times U}\E\big[|\Psi^{v,u}(t,b)-\Psi^{v,u}_T(t,b)|\big|\mcF_t\big]\to 0$, $\Prob$-a.s., as $T\to\infty$ for all $v\in\mcU^f$ uniformly in $u\in\mcU^k$.

By Lemma~\ref{lem:J-star-T} and uniform convergence it follows from Lemma~\ref{lem:mcH-prop}.(d) that $J^*:=\lim_{T\to\infty}J^*_T\in\mcH_\bbF$. The desired result now follows by Lemma~\ref{lem:mcH-prop}.(a) while noting that by the construction of $\ell$ in Assumption~\ref{ass:onSFDE}.(\ref{ass:onSFDE-c}), a simplified version of Lemma~\ref{lem:J-star-T} (without having to consider maximization over $u$) applied to each of the $\ell_i$ gives that there is an $h\in\mcH_\bbF$ such that $h(\tau,b)=-\E\big[\ell(\tau,X^{v}_\tau,b)|\mcF_\tau\big]$, $\Prob$-a.s.~(with an exception set that is independent of $b$).\qed\\

\begin{rem}
In a perfect information setting, \ie when $\bbF=\bbG$, we note that $\ell(t,x,b)$ can be taken to be any upper semi-continuous function in $b$ that satisfies the remaining properties of polynomial growth and local Lipschitz continuity.
\end{rem}

\subsection{The random horizon setting}
We turn instead to the reward
\begin{align}\label{ekv:rewardSFDE-RH}
J^{\eta}(u)=\E\bigg[\int_0^\eta e^{-\rho(t)}\phi(t,X^u_t)dt+e^{-\rho(\eta)}\psi(\eta,X^{[u]_{N(\eta-)}}_\eta)-\sum_{j=1}^Ne^{-\rho(\tau_j)}\ell(\tau_j,X^{\nu,j-1}_{\tau_j},\beta_j)\bigg].
\end{align}
where $\eta$ is a $\bbG$-stopping and $N(\eta-):=\sup\{j:\tau_j<\eta\}\vee 0$. A notable convention applied in \eqref{ekv:rewardSFDE-RH} is that the terminal reward disregards interventions made at the horizon. This is natural from an applications perspective as it is generally to late to intervene at a default in a financial setting or at the failure of a unit in an engineering application.

In addition to the requirements listed in Assumption~\ref{ass:onSFDE}, we make the following assumptions:
\begin{ass}\label{ass:onSFDE-RH}
The terminal reward $\psi:\R_+\times \R^d\to\R$ is Borel-measurable, satisfies the growth condition
\begin{align*}
  |\psi(t,x)|\leq C(1+|x|^q)
\end{align*}
and for each $L>0$ there is a $C>0$ such that
\begin{align*}
  |\psi(t,x)-\psi(t,x')|\leq C|x-x'|
\end{align*}
whenever $|x|\vee|x'|\leq L$. Moreover, if there is a sequence $(\theta_j)_{j\geq 0}$ in $\mcT^f$ such that $\theta_j\nearrow \eta$ on some set $B\in\mcG$, then there is a $\Prob$-null set $\mcN$ such that on $B\setminus\mcN$ we have for every $k$ and every $(y,\vecb)\in \bbD\times U^k$ that
\begin{align}
\psi(\eta,y_\eta)\geq \psi(\eta,\Gamma_{b_k}\circ\cdots\circ \Gamma_{b_1}(\eta,(y_s)_{s\leq \eta}))\label{ekv:psi-qlc}
\end{align}
where $\Gamma_{b}(\cdot,\cdot):=\Gamma(\cdot,\cdot,b)$ and $\circ$ denotes composition of functions.
\end{ass}

We introduce the following hypothesis:
\begin{hyp*}\textbf{[Disc.-B]} The terminal reward satisfies the following bound
\begin{align*}
\sup_{u\in\mcU^f}\E\Big[e^{-2\rho(\eta)}|\psi(\eta,X^u_\eta)|^2\Big]<\infty.
\end{align*}
Furthermore, for each $k\geq 0$ there is an $\epsilon>0$ such that for all $T\geq T'$ for some $T'>0$ we have
\begin{align*}
\sup_{u\in\mcU^k}\E\Big[\sup_{(t,b)\in \R_+\times U}\E\big[\ett_{[\eta\geq T]}e^{-\rho(\eta)}|\psi(\eta,X^{v\circ(t,b)\circ u}_{\eta})|\big|\mcF_t\big]\Big]\leq e^{-\epsilon T}
\end{align*}
for all $v\in\mcU^f$.
\end{hyp*}


We have the following extension of Proposition~\ref{prop:SFDE-OC}.

\begin{prop}\label{prop:SFDE-OC-RH}
Under Hypotheses~\textbf{Disc.-A} and \textbf{Disc.-B} there is a $u^*\in\mcU^f$ such that $J^{\eta}(u^*)\geq J^{\eta}(u)$ for all $u\in\mcU^f$. Furthermore, $u^*$ is given by the recursion \eqref{ekv:taujDEF}-\eqref{ekv:betajDEF}, with
\begin{align*}
\varphi(u):=\int_0^{\eta} e^{-\rho(t)}\phi(t,X^u_t)dt+e^{-\rho(\eta)}\psi(\eta,X^{[u]_{N(\eta-)}}_\eta)
\end{align*}
and
\begin{align*}
c(u\circ(t,b)):=e^{-\rho(t\vee\tau_N)}\ell(t\vee\tau_N,X^u_{t\vee\tau_N},b).
\end{align*}
If, in addition $\eta$ is an $\bbF$-stopping time, then $\tau^*_j<\eta$ for all $1\leq j\leq N^*$ on $\{\eta<\infty\}$.
\end{prop}

\noindent\emph{Proof.} We note that all details in the proof of Proposition~\ref{prop:SFDE-OC} transfer immediately to this situation except for the quasi-left upper semi-continuity property in the definition of $\mcH_\bbF$ (Definition~\ref{def:mcH}). When $\theta<\eta$ the assertion follows by Lemma~\ref{lem:SFDEflow} and the local Lipschitz property of $\psi$ and when $\theta>\eta$ the result is immediate. We thus assume that $\theta_j\nearrow\eta$ on some measurable set $B\subset\Omega$.
Then, we have
\begin{align*}
&\ett_{B}(\varphi(v\circ(\theta_j,b)\circ u)-\varphi(v\circ(\theta,b)\circ u))
\\
&\leq \int_{\theta_j}^{\theta} e^{-\rho(t)}|\phi(t,X^{v\circ(\theta_j,b)\circ u}_t)-\phi(t,X^{v\circ(\theta,b)\circ u}_t)|dt+\ett_{B}e^{-\rho(\eta)}(\psi(\eta,X_\eta^{v\circ(\theta_j,b)\circ u})-\psi(\eta,X_\eta^v)),
\end{align*}
where the first term on the right-hand side tends to zero, $\Prob$-a.s. Concerning the second term we have
\begin{align*}
\ett_{B}e^{-\rho(\eta)}(\psi(\eta,X_\eta^{v\circ(\theta_j,b)\circ u})-\psi(\eta,X_\eta^v))&\leq \ett_{B}e^{-\rho(\eta)}(\psi(\eta,\Gamma_{\beta_{N(\eta-)}}\circ\cdots\circ\Gamma_{\beta_1}\circ\Gamma_b(\eta,X_\eta^{v})-\psi(\eta,X_\eta^v))
\\
&\quad +e^{-\rho(\eta)}|\psi(\eta,X^{v\circ (\theta_j,b)\circ [u]_{N(\eta-)}}_\eta)-\psi(\eta,X^{v\circ (\theta,b)\circ [u]_{N(\eta-)}}_\eta)|,
\end{align*}
where the first term on the right-hand side is $\Prob$-a.s.~non-positive by Assumption~\ref{ass:onSFDE-RH} and the last term tends to zero, $\Prob$-a.s., by the local Lipschitz property of $\psi$ and Lemma~\ref{lem:SFDEflow} in combination with Proposition~\ref{prop:SFDEmoment}, the polynomial growth condition on $\psi$ and Hypothesis \textbf{Disc.-B}.

The last assertion follows by noting that since $c>0$ it will never be optimal to intervene at times greater than or equal to $\eta$.\qed\\

We note the following distinction between the finite (deterministic) horizon and the random horizon settings:

\begin{rem}
In the case when $\Prob(\eta=T)=1$ for some $T\geq 0$ it follows from the proof of Proposition~\ref{prop:SFDE-OC-RH} that we can relax \eqref{ekv:psi-qlc} to
\begin{align*}
\psi(T,y_T)\geq \psi(T,\Gamma_{b_k}\circ\cdots\circ \Gamma_{b_1}(T,(y_s)_{s\leq T}))-\sum_{j=1}^k \ell(T,\Gamma_{b_j}\circ\cdots\circ \Gamma_{b_1}(T,(y_s)_{s\leq T}),b_j).
\end{align*}
\end{rem}

To see that there is an actual distinction here consider the following example:
\begin{example}
We let $\mcF$ be the trivial $\sigma$-algebra $\{\emptyset,\Omega\}$ and assume that $\Prob(\eta=x)=\left\{\begin{array}{ll}0.5, & x=1 \\ 0.5, & x=2\end{array}\right.$. We take $U:=\{1\}$ and set $X_t=\ett_{[\tau_1,\infty)}(t)$. Then, with the rewards $\phi\equiv 0$, $\psi(t,x)=xe^{|t-1|}$, the intervention cost $\ell(t,x,b)=e^{|t-1|}$ and the discount $\rho\equiv 0$, we get
\begin{align*}
\sup_{u\in\mcU}J^\eta(u)=0.5(e^{1}-1),
\end{align*}
but there is no control that attains this value.
\end{example}


\bibliographystyle{plain}
\bibliography{ElepInfHor_ref}

\begin{thebibliography}{10}

\bibitem{AgramOksen}
N.~Agram and B.~{\O}ksendal.
\newblock Stochastic control of memory mean-field processes.
\newblock {\em Appl. Math. Optim.}, 79:181--204, 2019.

\bibitem{RAid2}
R.~A\"{i}d, S.~Federico, H.~Pham, and B.~Villeneuve.
\newblock Explicit investment rules with time-to-build and uncertainty.
\newblock {\em J. Econom. Dynam. Control}, 51:240--256, 2015.

\bibitem{BarIlan95}
A.~Bar-Ilan and A.~Sulem.
\newblock Explicit solution of inventory problems with delivery lags.
\newblock {\em Math. Oper. Res.}, 20(3), 1995.

\bibitem{BaseiImpulse}
M.~Basei.
\newblock Optimal price management in retail energy markets: an impulse control
  problem with asymptotic estimates.
\newblock {\em Math Meth Oper Res}, 89:355--383, 2019.

\bibitem{BensLionsImpulse}
A.~Bensoussan and J.L. Lions.
\newblock {\em Impulse Control and Quasivariational inequalities}.
\newblock Gauthier-Villars, Montrouge, France, 1984.

\bibitem{BertsekasShreve}
D.~P. Bertsekas and S.~E. Shreve.
\newblock {\em Stochastic optimal control: The discrete-time case}.
\newblock Academic Press, 1978.

\bibitem{Bruder}
B.~Bruder and H.~Pham.
\newblock Impulse control problem on finite horizon with execution delay.
\newblock {\em Stochastic Process. Appl.}, 2009.

\bibitem{CarmLud}
R.~Carmona and M.~Ludkovski.
\newblock Pricing asset scheduling flexibility using optimal switching.
\newblock {\em Appl. Math. Finance}, 15:405--447, 2008.

\bibitem{DelMeyer1}
C.~Dellacherie and P.-A. Meyer.
\newblock {\em Probabilit{\'e}s et Potentiel, I-IV}.
\newblock Hermann, Paris, 1975.

\bibitem{DelMeyer2}
C.~Dellacherie and P.-A. Meyer.
\newblock {\em Probabilit{\'e}s et Potentiel, V-VIII}.
\newblock Hermann, Paris, 1980.

\bibitem{DjehiceImpulse}
B.~Djehiche, S.~Hamad\'ene, and I.~Hdhiri.
\newblock Stochastic impulse control of non-markovian processes.
\newblock {\em Appl Math Optim}, 61(1):1--26, 2010.

\bibitem{DjehicheInfHorImp}
B.~Djehiche, S.~Hamad\'ene, I.~Hdhiri, and H.~Zaatra.
\newblock Infinite horizon stochastic impulse control with delay and random
  coefficients.
\newblock {\em arXiv:1904.11924}, 2019.

\bibitem{BollanMSwitch1}
B.~Djehiche, S.~Hamad\'ene, and A.~Popier.
\newblock A finite horizon optimal multiple switching problem.
\newblock {\em SIAM J. Control Optim.}, 47(4):2751--2770, 2009.

\bibitem{ElKarouiLN}
N.~El~Karoui.
\newblock {\em Les aspects probabilistes du contr{\^o}le stochastique}.
\newblock Ecole d'Et\'e de SaintFlour IX 1979. Lecture Notes in Math. Springer,
  Berlin., 1981.

\bibitem{ElKaroui1}
N.~El-Karoui, C.~Kapoudjian, E.~Pardoux, S.~Peng, and M.~C. Quenez.
\newblock Reflected solutions of backward {SDEs} and related obstacle problems
  for {PDEs}.
\newblock {\em Ann. Probab.}, 25(2):702--737, 1997.

\bibitem{HamRefBSDE}
S.~Hamad\'ene.
\newblock Reflected {BSDE's} with discontinuous barrier and application.
\newblock {\em Stochastics: An International Journal of Probability and
  Stochastic Processes}, 74(3-4):571--596, 2002.

\bibitem{HamZhang}
S.~Hamad\'ene and J.~Zhang.
\newblock Switching problem and related system of reflected backward {SDEs}.
\newblock {\em Stochastic Processes and their Applications}, 120(4):403--426,
  2010.

\bibitem{Hdhiri}
I.~Hdhiri and M.~Karouf.
\newblock Optimal stochastic impulse control with random coefficients and
  execution delay.
\newblock {\em Stochastics: An International Journal of Probability and
  Stochastic Processes}, 90(2):151--164, 2018.

\bibitem{JonteSFDE}
J.~J{\"o}nsson and M.~Perninge.
\newblock Finite horizon impulse control of stochastic functional differential
  equations.
\newblock {\em arXiv:2006.09768}, 2020.

\bibitem{ElKarouiTan}
N.~El Karoui and X.~Tan.
\newblock Capacities, measurable selection and dynamic programming part i:
  Abstract framework.
\newblock {\em arXiv:1310.3363}, 2013.

\bibitem{Korn}
R.~Korn.
\newblock Some applications of impulse control in mathematical finance.
\newblock {\em Math Meth Oper Res}, 50:493--518, 1999.

\bibitem{MartyrSigned}
R.~Martyr.
\newblock Finite-horizon optimal multiple switching with signed switching
  costs.
\newblock {\em Math. Oper. Res.}, 41(4):1432--1447, 2016.

\bibitem{OksenSulemBok}
B.~{\O}ksendal and A.~Sulem.
\newblock {\em Applied Stochastic Control of Jump Diffusions}.
\newblock Springer, 2007.

\bibitem{OksenImpulse}
B.~{\O}ksendal and A.~Sulem.
\newblock Optimal stochastic impulse control with delayed reaction.
\newblock {\em Appl. Math. Optim.}, 58:243--255, 2008.

\bibitem{PSImpulsive}
J.~Palczewski and L.~Stettner.
\newblock Impulsive control of portfolios.
\newblock {\em Appl Math Optim}, 56:67--103, 2007.

\bibitem{SwitchElephant}
M.~Perninge.
\newblock A finite horizon optimal switching problem with memory and
  application to controlled sddes.
\newblock {\em Math Meth Oper Res}, 2019.

\bibitem{Protter}
P.~Protter.
\newblock {\em Stochastic Integration and Differential Equations}.
\newblock Springer, Berlin, 2nd edition, 2004.

\end{thebibliography}
\end{document}